\documentclass{amsart}
\usepackage{graphicx}
\usepackage{amssymb}
\usepackage{amsmath}
\usepackage{amsthm}
\usepackage{subfigure}

\newtheorem{thm}{Theorem}[section]

\newtheorem{lem}[thm]{Lemma}
\newtheorem{prop}[thm]{Proposition}
\theoremstyle{definition}
\newtheorem{defn}[thm]{Definition}
\newtheorem{rem}[thm]{Remark}

\newtheorem*{defn*}{Definition}
\newtheorem*{rems*}{Remarks}
\newtheorem*{rem*}{Remark}

\numberwithin{equation}{section}

\DeclareMathOperator {\Symp} {Symp}

\DeclareMathOperator {\Span} {span}

\def \algrest {\left [\Symp (\mathbb R^{2n})\right ]_{N}}

\def \algrestall {\bigl [\Lambda ^2(\mathbb R^{2n})\bigr ]_N}

\def \algrestclosed {\bigl [ Z ^2(\mathbb R^{2n})\bigr ]_N}

\begin{document}

\title[Symplectic $T_7$, $T_8$ singularities] {Symplectic $T_7$, $T_8$ singularities and Lagrangian tangency orders}
\author{Wojciech Domitrz}
\address{Warsaw University of Technology\\
Faculty of Mathematics and Information Science\\
Plac Politechniki 1\\
00-661 Warsaw\\
Poland\\}

\email{domitrz@mini.pw.edu.pl}
\author{\.{Z}aneta Tr\c{e}bska}
\address{Warsaw University of Technology\\
Faculty of Mathematics and Information Science\\
Plac Politechniki 1\\
00-661 Warsaw\\
Poland\\}

\email{ztrebska@mini.pw.edu.pl}

\thanks{The work of WD was supported by Polish MNiSW grant no. N N201 397237
}



\maketitle

\begin{abstract}
We study the local symplectic algebra of curves. We use the method
of algebraic restrictions to classify symplectic $T_7$
singularities. We define discrete symplectic invariants - the
Lagrangian tangency orders and compare them with the index of isotropy. We use these invariants to distinguish symplectic singularities of classical  $T_7$ singularity. We also
give the geometric description of symplectic classes of the singularity.
\end{abstract}

\section{Introduction}

In this paper we study the symplectic classification of singular curves under the following equivalence:

\begin{defn} \label{symplecto}
Let $N_1, N_2$ be germs of subsets of symplectic space $(\mathbb{R}^{2n}, \omega)$. $N_1, N_2$ are \textbf{symplectically equivalent} if there exists a symplectomorphism-germ $\Phi:(\mathbb{R}^{2n}, \omega) \rightarrow(\mathbb{R}^{2n}, \omega)$ such that $\Phi(N_1)=N_2$.
\end{defn}

We recall that  $\omega$ is a symplectic form if $\omega$ is a
smooth nondegenerate closed 2-form, and $\Phi:\mathbb{R}^{2n}
\rightarrow\mathbb{R}^{2n}$ is a symplectomorphism if $\Phi$ is
diffeomorphism and $\Phi ^* \omega=\omega$.

\medskip

Symplectic classification of curves was first studied by V. I.
Arnold. In \cite{Ar2} V. I. Arnold discovered new symplectic
invariants of singular curves. He proved that the $A_{2k}$
singularity of a planar curve (the orbit with respect to standard
$\mathcal A$-equivalence of parameterized curves) split into
exactly $2k+1$ symplectic singularities (orbits with respect to
symplectic equivalence of parameterized curves). He distinguished
different symplectic singularities by different orders of tangency
of the parameterized curve to the \emph{nearest} smooth Lagrangian
submanifold. Arnold posed a problem of expressing these invariants
in terms of the local algebra's interaction with the symplectic
structure and he proposed to call this interaction the local
symplectic algebra.

In \cite{IJ1} G. Ishikawa and S. Janeczko classified symplectic singularities of curves in the $2$-dimensional symplectic space. All simple curves in this classification are quasi-homogeneous.

\medskip

We recall that a subset $N$ of $\mathbb R^m$ is \textbf{quasi-homogeneous} if there exists
a coordinate system $(x_1,\cdots,x_m)$ on $\mathbb R^m$ and
positive numbers $w_1,\cdots,w_m$ (called weights) such that for any
point $(y_1,\cdots,y_m)\in \mathbb R^m$ and any $t\in \mathbb R$
if $(y_1,\cdots,y_m)$ belongs to $N$ then a point
$(t^{w_1}y_1,\cdots,t^{w_m}y_m)$ belongs to $N$.

\medskip

A symplectic form on a $2$-dimensional manifold is a special case of a volume form on a smooth manifold. The generalization of results in \cite{IJ1} to volume-preserving classification of singular varieties and maps  in arbitrary dimensions was obtained in \cite{DR}. The orbit of action of all diffeomorphism-germs agrees with volume-preserving orbit or splits into two volume-preserving orbits (in the case $\mathbb K=\mathbb R$) for germs which satisfy a special weak form of quasi-homogeneity e.g. the weak quasi-homogeneity of varieties is a quasi-homogeneity with non-negative weights $w_i\ge0$ and $\sum_i w_i>0$.

Symplectic singularity is stably simple if it is simple and remains simple if the ambient symplectic space is symplectically embedded (i.e. as a symplectic submanifold) into a larger symplectic space. In \cite{K} P. A. Kolgushkin classified the stably simple symplectic singularities of parameterized curves (in the $\mathbb C$-analytic category). All stably simple symplectic singularities of curves are quasi-homogeneous too.

In \cite{DJZ2} new symplectic invariants of singular quasi-homogeneous  subsets of a symplectic space were explained by the algebraic restrictions of the symplectic form to these subsets.

\medskip

The algebraic restriction is an equivalence class of the following relation on the space of differential $k$-forms:

Differential $k$-forms $\omega_1$ and $\omega_2$ have the same
{\bf algebraic restriction} to a subset $N$ if
$\omega_1-\omega_2=\alpha+d\beta$, where $\alpha$ is a $k$-form
vanishing on $N$ and $\beta$ is a $(k-1)$-form vanishing on $N$.

\medskip

In \cite{DJZ2} the generalization of Darboux-Givental theorem (\cite{ArGi})
to germs of arbitrary subsets of the symplectic space was obtained. This result reduces
the problem of symplectic classification of germs of quasi-homo\-ge\-neous subsets to
the problem of classification of algebraic restrictions of symplectic
forms to these subsets. For non-quasi-homogeneous subsets there is one more cohomological invariant except the algebraic restriction (\cite{DJZ2}, \cite{DJZ1}). The dimension of the space of algebraic restrictions of closed $2$-forms to a $1$-dimensional quasi-homogeneous isolated complete
intersection singularity $C$ is equal to the multiplicity of $C$ (\cite{DJZ2}). In \cite{D} it was proved that the space of algebraic restrictions of closed $2$-forms to a $1$-dimensional (singular) analytic variety is finite-dimensional. In \cite{DJZ2} the method of algebraic restrictions was applied to various classification problems in a symplectic space. In particular the complete symplectic classification of classical $A-D-E$ singularities of planar curves and $S_5$ singularity were obtained. Most of different symplectic singularity classes were distinguished by new discrete symplectic invariants: the index of isotropy and the symplectic multiplicity.

\medskip

In this paper following ideas from \cite{Ar2} and \cite{D} we use new discrete symplectic invariants - the Lagrangian tangency orders (section \ref{tang-order}). Although this invariant has definition similar to the index of isotropy its nature is different. Since the Lagrangian tangency order takes into account the weights of quasi-homogeneity of curves it allows us to distinguish more symplectic classes in many cases. For example using the Lagrangian tangency order we are able to distinguish classes $E_6^3$ and $E_6^{4,\pm}$ of classical planar singularity $E_6$ which can not be distinguished nor by the isotropy index nor by the symplectic multiplicity. In the paper we also present other examples of singularities which can be distinguished only by the Lagrangian tangency order. On the other hand, there are  singularities which symplectic classes can be distinguished by the index of isotropy but not by the Lagrangian tangency order, for example  the parametric curve with semigroup $(3,7,11)$ and  $T_8$ singularity. These examples show that there are no simple relations between the Lagrangian tangency order and the index of isotropy even for the case of parametric curves.

\medskip\
We also obtain the complete symplectic classification of the classical isolated complete intersection singularity $T_7$ using the method of algebraic restrictions (Theorem \ref{T7-main}). We
calculate discrete symplectic invariants for this classification (Theorems \ref{lagr-t7}) 
and we present geometric descriptions of its symplectic orbits (Theorem  \ref{geom-cond-t7}).

The paper is organized as follows.  In Section \ref{discrete} we present known discrete symplectic invariants and introduce the Lagrangian tangency orders. We also compare the Lagrangian tangency order and the index of isotropy.   
Symplectic classification of $T_7$ singularity is studied in Section \ref{sec-t7}.
In Section \ref{proofs} we recall the method of algebraic restrictions and use it to classify $T_7$ symplectic singularities. 

\section{Discrete symplectic invariants.}\label{discrete}

We define discrete symplectic invariants to distinguish
symplectic singularity classes.  The first one is the symplectic
multiplicity (\cite{DJZ2}) introduced  in \cite{IJ1} as a
symplectic defect of a curve.

\medskip

Let $N$ be a germ of a subset of $(\mathbb R^{2n},\omega)$.

\begin{defn}
\label{def-mu}
 The {\bf symplectic multiplicity} $\mu_{sympl}(N)$ of  $N$ is the codimension of
 a symplectic orbit of $N$ in an orbit of $N$ with respect to the action of the group of local diffeomorphisms.
\end{defn}

The second one is the index of isotropy \cite{DJZ2}.

\begin{defn}
The {\bf index of isotropy} $ind(N)$ of $N$ is the maximal
order of vanishing of the $2$-forms $\omega \vert _{TM}$ over all
smooth submanifolds $M$ containing $N$.
\end{defn}

This invariant has geometrical interpretation.
An equivalent definition is as follows:
the index of isotropy of $N$ is the maximal order of tangency between
non-singular submanifolds containing $N$ and non-singular
isotropic submanifolds of the same dimension.  The index of
isotropy is equal to $0$ if $N$ is not contained in any
non-singular submanifold which is tangent to some isotropic
submanifold of the same dimension. If $N$ is contained in a
non-singular Lagrangian submanifold then the index of isotropy
is $\infty $.

\begin{rem} If $N$ consists of invariant components $C_i$ we can calculate index of isotropy for each component $ind(C_i)$ as the maximal
order of vanishing of the $2$-forms $\omega \vert _{TM}$ over all
smooth submanifolds $M$ containing $C_i$.
\end{rem}

The symplectic multiplicity and the index of isotropy can be described in terms of algebraic restrictions (Propositions \ref{sm} and \ref{ii} in Section \ref{proofs}).



\subsection{Lagrangian tangency order}
\label{tang-order}

There is one more discrete symplectic invariant introduced in \cite{D} following ideas from \cite{Ar2} which is defined specifically for a parameterized curve. This is the maximal
tangency order of a curve $f:\mathbb R\rightarrow M$ to a smooth Lagrangian submanifold. If $H_1=...=H_n=0$ define a smooth submanifold $L$ in the symplectic space then the tangency order of
a curve $f:\mathbb R\rightarrow M$ to $L$ is the minimum of the orders of vanishing at $0$ of functions $H_1\circ f,\cdots, H_n\circ f$. We denote the tangency order of $f$ to $L$ by $t(f,L)$.

\begin{defn}
The {\bf Lagrangian tangency order} $Lt(f)$\textbf{ of a curve} $f$ is the
maximum of $t(f,L)$ over all smooth Lagrangian submanifolds $L$ of
the symplectic space.
\end{defn}

The Lagrangian tangency order of a quasi-homogeneous curve in a symplectic space can also be  expressed in terms of algebraic restrictions  (Proposition \ref{lto} in Section \ref{proofs}).


We can generalize this invariant for curves which may be
parameterized analytically. Lagrangian tangency order is the same
for every 'good' analytic parameterization of a curve \cite{W}.
Considering only such parameterizations we can choose one and
calculate the invariant for it. It is easy to show that this
invariant doesn't depend on chosen parameterization.

\begin{prop}
Let $f:\mathbb R\rightarrow M$ and $g:\mathbb R\rightarrow M$ be good analytic parameterizations of the same curve. Then $Lt(f)=Lt(g)$.
\end{prop}

\begin{proof}

There  exists a diffeomorphism $\theta: \mathbb R\rightarrow \mathbb R$ such that $g(s)=f(\theta(s))$ and $\frac{d\theta}{ds}|_0\ne 0$. Let $H_1=\ldots=H_n=0$ define a smooth submanifold $L$ in the symplectic space. If $\frac{d^l(H_i\circ f)}{dt^l}|_0=0$ for $l=1,...,k$ then $$\frac{d^{k+1}(H_i\circ g)}{ds^{k+1}}|_0=\frac{d^{k+1}(H_i\circ f\circ \theta)}{ds^{k+1}}|_0=\frac{d^{k+1}(H_i\circ f)}{dt^{k+1}}|_0\cdot (\frac{d\theta}{ds})^{k+1}|_0$$ so  the orders of vanishing at $0$ of functions $H_i\circ f$ and $H_i\circ g$ are equal and hence $t(f,L)=t(g,L)$, which implies that $Lt(f)=Lt(g)$.
\end{proof}

\medskip

We can generalize Lagrangian tangency order for sets containing parametric curves.
Let $N$ be a subset of a symplectic space $(\mathbb R^{2n},\omega)$.

\begin{defn}
The {\bf  tangency order of the germ of a subset $N$ to the germ
of a submanifold $L$} $t[N,L]$ is equal to the minimum of $t(f,L)$
over all parameterized curve-germs $f$ such that $Im f\subseteq
N$.
\end{defn}

\begin{defn}
The {\bf Lagrangian tangency order of $N$} $Lt(N)$ is equal to the
maximum of $t[N,L]$ over all smooth Lagrangian submanifold-germs
$L$ of the symplectic space.
\end{defn}

In this paper we consider $N$ which are singular analytic curves.
They may be identified with a multi-germ of parametric curves. We
define invariants which are special cases of the above definition.

\medskip

\noindent Consider a multi-germ $(f_i)_{i\in\{1,\cdots,r\}}$ of analytically parameterized curves $f_i$. For any smooth  submanifold $L$ in the symplectic space we have $r$-tuples $(t(f_1,L), \cdots, t(f_r,L))$.

\begin{defn}
For any $I\subseteq \{1,\cdots, r\}$ we define \textbf{the tangency order of the multi-germ } $(f_i)_{i\in I}$ to $L$:
$$t[(f_i)_{i\in\ I},L]=\min_{i\in\ I} t(f_i,L).$$
\end{defn}

\begin{defn}
The {\bf Lagrangian tangency order} $Lt((f_i)_{i\in\ I})$ \textbf{of a multi-germ } $(f_i)_{i\in I}$ is the maximum of $t[(f_i)_{i\in\ I},L]$ over all smooth Lagrangian submanifolds $L$ of the symplectic space.
\end{defn}

\medskip

For multi-germs we can also define relative invariants according to selected branches or collections of branches.

\begin{defn}
 Let $S\subseteq I\subseteq \{1,\cdots,r\}$. For $i\in S$ let us fix numbers $t_i\leq Lt(f_i)$. The {\bf relative Lagrangian tangency order} $Lt[(f_i)_{i\in I}:(S,(t_i)_{i\in S})]$ of a multi-germ  $(f_i)_{i\in I}$  related to $S$ and $(t_i)_{i\in S}$ is the maximum of $t[(f_i)_{i\in I\setminus S},L]$ over all smooth Lagrangian submanifolds $L$ of the symplectic space for which $t(f_i,L)=t_i$, if such submanifolds exist, or $-\infty$ if there are no such submanifolds.
\end{defn}

We can also define special relative invariants according to selected branches of multi-germ.

\begin{defn}
For fixed $j\in I$  the \textbf{ Lagrangian tangency order related
to}  $f_j$ of a multi-germ  $(f_i)_{i\in I}$ \textbf{} denoted by
$Lt[(f_i)_{i\in I}: f_j]$ is the maximum of $t[(f_i)_{i\in
I\setminus\{j\}},L]$ over all smooth Lagrangian submanifolds $L$
of the symplectic space for which $t(f_j,L)=Lt(f_j)$,
\end{defn}

These invariants have geometric interpretations. If
$Lt(f_i)=\infty$ then a branch $f_i$ is included in a smooth
Lagrangian submanifold. If $Lt((f_i)_{i\in\ I})=\infty$ then
exists a Lagrangian submanifold containing all curves $f_i$ for
$i\in I$.

We may use these invariants to distinguish symplectic singularities.

\subsection{Comparison of the Lagrangian tangency order and the index of isotropy}
\label{ind-order}
Definitions of the Lagrangian tangency order and the  index of isotropy are similar. They show how far a variety $N$ is from the nearest non-singular Lagrangian submanifold. The index of isotropy of a quasi-homogeneous set $N$ is $\infty$ if and only if the Lagrangian tangency order of $N$ is $\infty$. Studying classical singularities we have found examples of all possible interactions between these invariants.

\bigskip

$\bullet$ \ \ For some singularities the index of isotropy distinguishes the same symplectic classes which can be distinguished by the Lagrangian tangency order. It is observed for example for  planar curves - the classical $A_k$ and $D_k$ singularities (Tables  \ref{tab-a} and \ref{tab-d}) and for $S_{\mu}$ singularities studied in \cite{DT}.

\medskip

 A complete symplectic classification of classical $A-D-E$ singularities of planar curves was obtained using a method of algebraic restriction in \cite{DJZ2}. Below we compare the Lagrangian tangency order and the index of isotropy for these singularities. A curve $N$ may be described as a parameterized curve or as a union of parameterized components $C_i$ preserved by local diffeomorphisms in the symplectic space $(\mathbb R^{2n}, \omega _0=\sum_{i=1}^ndp_i\wedge dq_i)$, $n\ge 2$. For calculating the Lagrangian tangency orders we  give their parameterization in the coordinate system $(p_1,q_1,p_2,q_2,\cdots,p_n,q_n)$.

 \medskip

 Denote by $(A_k)$ the class of varieties in  a fixed symplectic space $(\mathbb R^{2n}, \omega )$ which are diffeomorphic to
\begin{equation}
\label{defAk} A_k=\{x\in \mathbb R ^{2n\geq 4}\,:x_1^{k+1}-x_2^2=x_{\geq 3}=0\}.\end{equation}
  A curve $N\in(A_k)$ can be described as parameterized singular curve $C$ for $k$ - even, or as a pair of two smooth parameterized branches $B_+$ and $B_-$ if $k$ is odd. By $Lt(N)$ we denote $Lt(C)$ or $Lt(B_+,B_-)$ respectively.

\renewcommand*{\arraystretch}{1.3}
\begin{small}
\begin{table}[h!]
\begin{tabular}{|p{3.0cm}|p{5.0cm}|p{1.3cm}|p{0.8cm}|}

\hline

Normal form & Parameterization & $Lt(N)$ & $ind$ \\ \hline

$A_k^{0\le i\le k-1}$ ($k$ - even)
  &  $C:(t^2,t^{k+1+2i}, t^{k+1},0,\!\cdots\!,0)$ & $k+1+2i$ & $i$  \\ \hline
$A_k^k $ \;\;\; ($k$ - even)& $C:(t^2,0, t^{k+1},0,\!\cdots\!,0)$ & $\infty $ &  $\infty $ \\ \hline

$A_k^{0\le i\le k-1}$  \ ($k$ - odd)
 & $B_{\pm}\!:(t, \pm t^{\frac{k+1}{2}+i},\pm t^{\frac{k+1}{2}},0,\cdots\!,0)$ & $\frac{k+1}{2}+i$ & $i$  \\ \hline

$A_k^k $, \;\;\;\; ($k$ - odd) & $B_{\pm}\!:(t, 0,\pm t^{\frac{k+1}{2}},0,\cdots\!,0)$  & $\infty $ &  $\infty $  \\ \hline

\end{tabular}

\medskip

\caption{\small Comparison of  symplectic invariants of $A_k$ singularity.}\label{tab-a}
\end{table}
\end{small}






 \medskip
 \newpage

 Denote by $(D_k)$ for $k\geq 4$ the class of varieties in  a fixed symplectic space $(\mathbb R^{2n}, \omega )$ which are diffeomorphic to
\begin{equation}
\label{defDk} D_k=\{x\in \mathbb R ^{2n\geq 4}\,:x_1^2x_2 - x_2^{k-1}=x_{\geq 3}=0\}.\end{equation}
 A curve $N\in(D_k)$ consists of $2$ invariant components: $C_1$ - smooth and $C_2$ - singular diffeomorphic to $A_{k-3}$. $C_2$ may consists of one or two branches depending on $k$. To distinguish the symplectic classes completely we need two invariants: $Lt(N)$ -- the Lagrangian tangency order of $N$ and $Lt(C_2)$ -- the Lagrangian tangency order of the singular component $C_2$. Equivalently we can use the index of isotropy of $N$ -- $ind$ and the index of isotropy of $C_2$ -- $ind_2$. 
\renewcommand*{\arraystretch}{1.3}
\begin{small}
\begin{table}[h]

\begin{tabular}{|p{1.3cm}|p{4.0cm}|l|l|l|l|}

\hline

Normal \newline form & $f(t)$ & $Lt(N)$ &  $Lt(C_2)$  & $ind$ & $ind_2$ \\ \hline

$D_k^0$ &  $ t^{2 \lambda_k}$ & $2 \lambda_k$ & $(k-2) \lambda_k$ & $0$ & $0$ \\ \hline

$D_k^1$   & $ bt^{k \lambda_k} +\frac{1}{2}t^{4 \lambda_k}$ & $k \lambda_k$  & $k \lambda_k$ & $1$ & $1$ \\ \hline

$D_k^i$ & $ bt^{k \lambda_k}+\frac{1}{i+1}t^{2(i+1)\lambda_k}$, $b\!\neq\! 0$ & $ k \lambda_k$ & $ (k\!-\!2\!+\!2i)\lambda_k$ & $1$ & $i$ \\ \cline{2-6}
 $^{1< i< k-3}$  & $ \frac{1}{i+1}t^{2(i+1)\lambda_k}$  &  $(k\!-\!2\!+\!2i)\lambda_k$ & $(k\!-\!2\!+\!2i)\lambda_k$ & $i$ & $i$ \\ \hline

$D_k^{k-3,\pm} $  & $ (\pm 1)^k t^{k \lambda_k}+\!\frac{b}{k-2}t^{2(k-2)} \lambda_k$ & $k \lambda_k$ & $\infty $ &  $1$ & $\infty$ \\ \hline

 $D_k^{k-2} $  & $ \frac{1}{k-2}t^{2(k-2) \lambda_k}$ & $ (3k-8)\lambda_k$ & $\infty $ &  $k-3$ & $\infty$ \\ \hline

 $D_k^{k-1} $  & $ \frac{1}{k-1}t^{2(k-1) \lambda_k}$ & $ (3k-6)\lambda_k$ & $\infty $ &  $k-2$ & $\infty$ \\ \hline

 $D_k^k $  & $0$ & $\infty $ &  $\infty $ & $\infty $ & $\infty$ \\ \hline

\end{tabular}

\medskip

\caption{\small Symplectic invariants of $D_k$ singularity. 
The branch  $C_1$  has a form
$(t,0,0,0,\cdots,0)$. If $k$ is odd then $C_2$ has a form $(t^{k-2},f(t),t^2, 0,\cdots,0)$ and $\lambda_k=1$. If $k$ is even then $C_2$ consists of two branches: $B_{\pm}: (\pm t^{(k-2)/
2},f(t),t,0,\cdots,0)$  and $\lambda_k=\frac{1}{2}$.}\label{tab-d}
\end{table}
\end{small}

\medskip


 $\bullet$ \ \ There are also symplectic singularities distinguished by the Lagrangian tangency order but not by the index of isotropy. The simplest example is planar  singularity $E_6$ (Table \ref{tab-e6}). Such a "more sensitivity" of the Lagrangian tangency order we also observe for $E_7$ and $E_8$ singularities and for parametric curves with semigroups $(3,4,5)$, $(3,5,7)$ and $(3,7,8)$ studied in \cite{D}. 

\medskip

Denote by $(E_6)$  the class of varieties in  a fixed symplectic space $(\mathbb R^{2n}, \omega )$ which are diffeomorphic to
\begin{equation}
\label{defE6} E_6=\{x\in \mathbb R ^{2n\geq 4}\,:x_1^3 - x_2^4=x_{\geq 3}=0\}.\end{equation}

\renewcommand*{\arraystretch}{1.2}
\begin{table}[h!bt]
\begin{tabular}{|l|l|l|l|l|}\hline

Normal form & Parameterization & $Lt(N)$ & $ind$ & $\mu ^{\rm symp}$\\ \hline

$E_6^0$ & $(t^4, t^3,t^3,0,\cdots,0)$& $4$ & 0 & 0 \\ \hline

$E_6^{1,\pm}$ &  $(t^4, \pm \frac{1}{2}t^6+bt^7,t^3,0,\cdots,0)$&  $7$  & $1$ & $2$  \\ \hline

$E_6^2$ &  $(t^4, t^7+\frac{b}{3}t^9,t^3,0,\cdots,0)$& $8$  & $1$ & $3$  \\ \hline

$E_6^3$ &  $(t^4, \frac{1}{3}t^9+\frac{b}{2}t^{10}, t^3,0,\cdots,0)$ & $\mathbf{10}$ & $\mathbf{2}$ & $\mathbf{4}$  \\ \hline

$E_6^{4,\pm}$ &  $(t^4, \pm \frac{1}{2}t^{10},t^3,0,\cdots,0)$& $\mathbf{11}$ & $\mathbf{2}$ & $\mathbf{4}$  \\ \hline

$E_6^5$ &  $(t^4, \frac{1}{3}t^{13},t^3,0,\cdots,0)$& $14$ & $3$ & $5$  \\ \hline

$E_6^6$ &  $(t^4,0,t^3,0,\cdots,0)$& $\infty$ & $\infty $  & $6$ \\ \hline 
\end{tabular}

\medskip
\caption{\small Symplectic invariants of $E_6$ singularity.}\label{tab-e6}
\end{table}

\medskip

 As we see in  Table \ref{tab-e6} we are able to distinguish by the Lagrangian tangency order the classes $E_6^3$ and $E_6^{4,\pm}$ which can not be distinguished nor by the index of isotropy nor by the symplectic multiplicity.

\bigskip


$\bullet$ \ \ Some symplectic singularities can be distinguished by the index of isotropy but not by the Lagrangian tangency order. Such situation we observe for a parametric quasi-homogeneous curve-germ  with semigroup $(3,7,11)$ listed as a stably simple singularity of curves in \cite{Ar1}. Another example is  $T_8$ singularity presented below (Table \ref{tabt8-lagr2} rows for $(T_8)^4$ and $(T_8)^{6,2}$).

\bigskip

The germ of a curve $f:(\mathbb{R},0)\rightarrow (\mathbb{R}^{2n},0)$ with semigroup $(3,7,11)$ is diffeomorphic to the curve $t\rightarrow (t^3,t^7, t^{11},0,\ldots,0)$. Among symplectic singularities of this curve-germ  in the symplectic space  $(\mathbb R^{2n},\omega=\sum_{i=1}^n dp_i \wedge dq_i)$ with the canonical coordinates $(p_1,q_1,\ldots,p_n,q_n)$ we have for example the  classes represented by the following normal forms:

\medskip

\renewcommand*{\arraystretch}{1.2}
\begin{center}
\begin{tabular}{cp{7cm}cc} 

Class & Normal form of $f$ & $Lt(f)$ & $ind$ \\ \hline

(1) & $t\rightarrow (t^3, t^{10},t^7,0,t^{11},0,\cdots,0)$& $10$ & 1  \\ 

(2) & $t\rightarrow (t^3, t^{11},t^7,0,t^{11},0,\cdots,0)$& $11$ & 0  \\ 

(3) &  $t\rightarrow (t^3, t^{10}+ct^{11},t^7,0,t^{11},0,\cdots,0), \ c\ne 0$&  $10$  & $0$   \\ \hline

\end{tabular}
\end{center}

\bigskip

Symplectic classes (1) and (3) have the same Lagrangian tangency order equal to $10$ but have different indices of isotropy -- $1$ and $0$ respectively. Symplectic classes (2) and (3) have the same index of isotropy equal to $0$ but have different Lagrangian tangency orders -- $11$ and $10$ respectively.  We also observe that the Lagrangian tangency order for class (1) is less than for class (2) but the inverse inequality is satisfied for the indices of isotropy.

\bigskip

Another example is $T_8$ singularity. Denote by $(T_8)$  the class of varieties in  a fixed symplectic space $(\mathbb R^{2n}, \omega )$ which are diffeomorphic to
\begin{equation}
\label{deft8} T_8=\{x\in \mathbb R ^{2n\geq 4}\,:x_1^2+x_2^3-x_3^4=x_2 x_3=x_{\geq 4}=0\}.\end{equation}

This is the classical $1$-dimensional isolated complete
intersection singularity $T_8$ (\cite{G}, \cite{AVG}).

Let $N\!\in\!(T_8)$. $N$ is quasi-homogeneous with weights $w(x_1)\!=6,  w(x_2)\!=4$, \ $w(x_3)\!=3$.
 A curve $N$ consists of $2$ invariant singular components: $C_1$ -- diffeomorphic to $A_2$ singularity and $C_2$ -- diffeomorphic to $A_3$ singularity which is a union of  two smooth branches $B_{+}$ and $B_{-}$. In local coordinates they have the form
\[\mathcal{C}_{1}=\{ x_1^2+x_2^3=0,\; x_3=x_{\geq 4}=0\},\]
\[\mathcal{B}_{\pm}=\{ x_1\pm x_3^2=0,\; x_2=x_{\geq 4}=0\}.\]

\smallskip

Using the method of algebraic restrictions  one can  obtain in the same way as it is presented in last two sections for the case of $T_7$ singularity the following complete classification of symplectic $T_8$ singularities.

\begin{thm}\label{t8-main}
Any stratified submanifold of the symplectic space  $(\mathbb R^{2n},\omega=\sum_{i=1}^n dp_i \wedge dq_i)$ which is diffeomorphic to $T_8$ is symplectically equivalent to one and only one of the
normal forms $(T_8)^i, i = 0,1,\cdots ,8$
. The parameters $c, c_1, c_2, c_3$ of the normal forms are moduli.

\medskip

\noindent $T_8$$^0: p_1^2 + p_2^3 - q_1^4 = 0, \  p_2q_1 = 0, \ q_2 = c_1q_1 - c_2p_1 \ p_{\ge 3} = q_{\ge 3} = 0, \  c_1\cdot c_2\ne 0$;

\medskip

\noindent $T_8$$^1_2: p_1^2 + p_2^3 - q_1^4 = 0,   p_2q_1 = 0,
q_2 = c_1q_1 - c_2 p_1 -c_3 p_1p_2,   p_{\ge 3} = q_{\ge 3} = 0,   c_1\cdot c_2=0$;

\medskip

\noindent $T_8$$^1_3: p_1^2 + q_1^3 - q_2^4 = 0, \  q_1q_2 = 0, \
p_2 = c_1q_1 + c_2p_1q_2, \  p_{\ge 3} = q_{\ge 3} = 0, \  c_1\cdot c_2\ne 0$;

\medskip

\noindent $T_8$$^2_3: p_1^2 + q_1^3 - q_2^4 = 0, \  q_1q_2 = 0, \
p_2 = c_1q_1+c_2p_1q_2+c_3p_1q_2^2, \  p_{\ge 3} = q_{\ge 3} = 0,$\par \ \ \ \ \ \  $c_1\cdot c_2= 0$;

\medskip

\noindent $T_8$$^2_{>3}: p_2^2 + p_1^3 - q_1^4 = 0, \  p_1q_1 = 0, \  q_2 = \frac{c_1}{2}q_1^2+\frac{c_2}{2}p_1^2,
 \   p_{\ge 3} = q_{\ge 3} = 0, \  c_1 \ne 0;$

\medskip

\noindent $T_8$$^{3,0}: p_2^2 + p_1^3 - q_1^4 = 0,    p_1q_1 = 0,   q_2 = \frac{c_1}{2}p_1^2+\frac{c_2}{3}q_1^3,
    p_{\ge 3} = q_{\ge 3} = 0,   (c_1,c_2) \ne (0,0);$

\medskip

\noindent $T_8$$^{5,0}:  p_2^2 + p_1^3 - q_1^4 = 0, \   p_1q_1 = 0, \  q_2 = \frac{c}{4}q_1^4,
 \   p_{\ge 3} = q_{\ge 3} = 0;$

\medskip

\noindent $T_8$$^{3,1}:  p_1^2 + p_2^3 - p_3^4 = 0, \   p_2p_3 = 0, \  q_1 = \frac{1}{2}p_3^2+\frac{c_2}{2}p_2^2,
  \ q_2=-c_1 p_1 p_3, \ p_{\ge 4} = q_{\ge 3} = 0;$

  \medskip

\noindent $T_8$$^4 \ : \ p_1^2 + p_2^3 - p_3^4 = 0, \   p_2p_3 = 0, \  q_1 = \frac{c_1}{2}p_2^2+\frac{c_2}{3}p_3^3,
 \  q_2=- p_1 p_3, \  p_{\ge 4} = q_{\ge 3} = 0,$\par \ \ \ \ \ $(c_1,c_2)\ne (0,0);$

  \medskip

\noindent $T_8$$^{6,1}:  p_1^2 + p_2^3 - p_3^4 = 0, \   p_2p_3 = 0, \  q_1 = \frac{c}{4}p_3^4,
 \  q_2=- p_1 p_3, \  p_{\ge 4} = q_{\ge 3} = 0;$

  \medskip

\noindent $T_8$$^{5,1}:  p_1^2 + p_2^3 - p_3^4 = 0, \   p_2p_3 = 0, \  q_1 = \frac{1}{2}p_2^2+\frac{c}{3}p_3^3,
   \ p_{\ge 4} = q_{\ge 2} = 0;$

  \medskip

\noindent $T_8$$^{6,2} :  p_1^2 + p_2^3 - p_3^4 = 0, \   p_2p_3 = 0,   \ q_1 = \frac{1}{3}p_3^3,
   \ p_{\ge 4} = q_{\ge 2} = 0;$

\medskip

\noindent $T_8$$^{7} \ : \  p_1^2 + p_2^3 - p_3^4 = 0, \   p_2p_3 = 0,  \ q_1 = \frac{1}{4}p_3^4,
   \ p_{\ge 4} = q_{\ge 2} = 0;$

\medskip

\noindent $T_8$$^{8} \ : \ p_1^2 + p_2^3 - p_3^4 = 0, \  p_2p_3 = 0,  \ q_{\ge 1} = p_{\ge 4} = 0.$

\end{thm}

Lagrangian tangency orders and indices of isotropy were used to obtain a  detailed  classification of $(T_8)$. A curve $N\in (T_8)$ may be described as a union of three parametrical branches $C_1, B_{+}, B_{-}$. Their parameterization  in the coordinate system $(p_1,q_1,p_2,q_2,\cdots,p_n,q_n)$ is presented in the second column of Tables \ref{tabt8-lagr1} and \ref{tabt8-lagr2}. To distinguish the classes of this singularity we need the following three invariants:
\begin{itemize}
  \item $Lt(N)=Lt(C_1,B_+,B_-)=\max\limits _L (\min \{t(C_1,L),t(B_+,L),t(B_-,L)\})$;
  \item $L_1=Lt(C_1)=\max\limits _L (t(C_1,L))$;
  \item $L_2=Lt(C_2)=\max\limits _L (\min \{t(B_+,L),t(B_-,L)\})$
\end{itemize}

where $L$ is a smooth Lagrangian submanifold of the symplectic space.

Branches $B_{+}$ and $B_{-}$ are diffeomorphic and are not preserved by all symmetries of $T_8$ so we can not use neither  $Lt(B_{+})$ nor $Lt(B_{-})$ as invariants.
 Considering the triples $(Lt, L_1, L_2)$ we obtain more detailed classification of symplectic singularities of $T_8$ than the classification given in Theorem \ref{t8-main}. Some subclasses appear having a natural geometric interpretation.

 We calculate  also index of isotropy of $N\in (T_8)$ denoted by $ind$ and the indices of isotropy of components $C_1$ and $C_2$ denoted respectively by $ind_1$ and $ind_2$. In Tables \ref{tabt8-lagr1} and \ref{tabt8-lagr2} we present  the comparison of the invariants.

\renewcommand*{\arraystretch}{1.3}
\begin{center}
\begin{table}[h]

    \begin{footnotesize}
    \noindent
    \begin{tabular}{|l|p{4cm}|l|c|c|c|c|c|c|}
                     \hline
    Class &  Parameterization  & Conditions & $Lt$   &  $L_1$ & $L_2$ & $ind$ & $ind_1$ & $ind_2$ \\ \hline
  $(T_8)^0$   &  $(t^3,0,-t^2,-c_2t^3,0,\cdots )$
  \newline $(\pm t^2,t,0,c_1 t\mp c_2t^2,0,\cdots )$ &  $c_1\cdot c_2\ne 0$& $2$ & $3$& $2$  &0 & 0&0 \\  \hline

  $(T_8)^1_2$ & $(t^3,0,-t^2,-c_2t^3+c_3t^5,0,\cdots )$ &   $c_1=0, c_2\ne 0 $ & $2$ & $3$ & $2$&0 &0 &0 \\ \cline{3-9}
  & $(\pm t^2,t,0,c_1 t\mp c_2t^2,0,\cdots )$ & $c_2=0, c_3\ne 0$  & $2$ & $5$ & $2$ &0 &1 &0\\ \cline{3-9}
         & & $c_2=c_3=0$ & $2$ & $\infty$ & $2$ &0 &$\infty$ &0 \\

                                         \hline

$(T_8)^1_3$  & $(t^3,-t^2,-c_1t^2, 0,0,\cdots)$ &  $c_1\cdot c_2\ne 0$ & $2$ & $3$ & $3$  &0 &0 &1 \\
      & $(\pm t^2,0,\pm c_2t^3,t,0,\cdots )$ &  &  &  & & & & \\ \hline

    $(T_8)^2_3$ &  $(t^3, -t^2, -c_1t^2,0, 0,\cdots)$ &  $c_1=0, c_2 \ne 0$ & $2$ & $3$ & $3$ &0 & 0&1 \\ \cline {3-9}
     & $(\pm t^2, 0, \pm c_2t^3\pm c_3t^4, t, 0, \cdots)$ & $ c_2=0, c_3\ne 0$ & $2$ & $3$ & $4$ &0 &0 &2 \\
     \cline{3-9}
         &  & $c_2=0, c_3= 0$ & $2$ & $3$ & $\infty$ &0 &0 &$\infty$ \\
          \hline

    $(T_8)^2_{>3}$ & $(-t^2,0,t^3,\frac{c_2}{2}t^4,0,\cdots )$ & $c_1\cdot c_2\ne 0$  & $2$ & $5$ & $3$&0 &1 &1 \\ \cline{3-9}
         & $(0,t,\pm t^2,\frac{c_1}{2}t^2,0,\cdots )$ & $c_1\ne 0, c_2=0$ & $2$ & $\infty$ & $3$&0 & $\infty$ & 1  \\ \hline

    $(T_8)^{3,0}$& $(-t^2,0,t^3,\frac{c_1}{2}t^4,0,\cdots )$ & $c_1\cdot c_2\ne 0$  & $2$ & $5$ & $4$ &0 &1 &2\\ \cline{3-9}
         & $(0,t,\pm t^2,\frac{c_2}{3}t^3,0,\cdots )$ & $c_1\ne 0, c_2=0$ & $2$ & $5$ & $\infty$ & 0&1 & $\infty$ \\ \cline{3-9}
         &  & $c_1=0, c_2\ne 0$  & $2$ & $\infty$ & $4$ &0 &$\infty$ &2\\ \hline

     $(T_8)^{5,0}$& $(-t^2,0,t^3,0,0,\cdots )$ &  & $2$ & $\infty$ & $\infty$&0 &$\infty$ &$\infty$ \\
        & $(0,t,\pm t^2,\frac{c}{4}t^4,0,\cdots )$ &   & & & & & &\\ \hline

\end{tabular}

\smallskip

\caption{\small Symplectic  invariants for symplectic classes of $T_8$ singularity when $\omega\vert_ W \ne 0$; $W$ - the tangent space to a non-singular $3$-dimensional manifold  in $(\mathbb{R}^{2n\ge 4},\omega)$ containing $N\in (T_8)$.}\label{tabt8-lagr1}

\end{footnotesize}
\end{table}
\end{center}

\begin{rem} We can notice that considering the pairs $(L_1,L_2)$ gives the same classification as considering the pairs (ind$_1$,ind$_2$).
To distinguish classes $(T_8)^0$ and $(T_8)^1_2$ for $c_2\ne 0, c_1=0$ we may use Lagrangian tangency order related to component $C_1$. We have $Lt[C_2:C_1]=1$ for class $(T_8)^0$ but $Lt[C_2:C_1]=2$ for class $(T_8)^1_2$ if $c_2\ne 0, c_1=0$. In similar way we can distinguish classes $(T_8)^1_3$ and $(T_8)^2_3$ for $c_2\ne 0, c_1=0$.
\end{rem}

\renewcommand*{\arraystretch}{1.3}
\begin{center}
\begin{table}[h]

    \begin{footnotesize}
    \noindent
    \begin{tabular}{|l|p{3.8cm}|l|c|c|c|c|c|c|}
                     \hline
    Class &  Parameterization  & Conditions &$Lt$   &  $L_1$ & $L_2$ & $ind$ & $ind_1$ & $ind_2$ \\ \hline
    $(T_8)^{3,1}$ & $(t^3,\frac{c_2}{2}t^4,-t^2,0,0,0,\cdots )$ & $ c_2\ne 0$ & $3$ & $5$ & $3$ &1 &1 &1\\ \cline{3-9}
     & $(\pm t^2,\frac{1}{2}t^2,0,\mp c_1t^3,t,0,\cdots )$ & $ c_2=0$  & $3$ & $\infty$ & $3$ &1 &$\infty$ &1 \\ \hline

      $(T_8)^4$  & $(t^3,\frac{c_1}{2}t^4,-t^2,0,0,0,\cdots )$ &   $c_1\cdot c_2\ne 0$ & $4$ & $5$ & $4$ &1 &1 &2 \\ \cline{3-9}
      & $(\pm t^2,\frac{c_2}{3}t^3,0,\mp t^3,t,0,\cdots )$ & $c_1=0, c_2\ne 0$ & $\mathbf{4}$ & $\mathbf{\infty}$ & $\mathbf{4}$ &$\mathbf{1}$ &$\mathbf{\infty}$ & $\mathbf{2}$\\ \cline{3-9}
      &  & $c_1\ne 0, c_2=0$  & $5$ & $5$ & $\infty$ &1 & 1&$\infty$ \\ \hline

      $(T_8)^{6,1}$  & $(t^3,0,-t^2,0,0,0,\cdots )$  &  & $5$ & $\infty$ & $\infty$ &1 &$\infty$ &$\infty$ \\
          & $(\pm t^2,\frac{c}{4}t^4,0,\mp t^3,t,0,\cdots )$ &   &  &  &  & & & \\ \hline

    $(T_8)^{5,1}$ & $(t^3,\frac{1}{2}t^4,-t^2,0,0,0,\cdots )$ & $c\ne 0$  & $4$ & $5$ & $4$ &1 &1 &2 \\ \cline{3-9}
       & $(\pm t^2,\frac{c}{3}t^3,0,0,t,0,\cdots )$ & $c=0$ & $5$ & $5$ & $\infty$ &1 &1 &$\infty$  \\ \hline

    $(T_8)^{6,2}$ & $(t^3,0,-t^2,0,0,0,\cdots )$ &   & $\mathbf{4}$ & $\mathbf{\infty}$ & $\mathbf{4}$ &$\mathbf{2}$ &$\mathbf{\infty}$ &$\mathbf{2}$ \\
         & $(\pm t^2,\frac{1}{3}t^3,0,0,t,0,\cdots )$ & &  &  & & & &   \\ \hline

     $(T_8)^{7}$ & $(t^3,0,-t^2,0,0,0,\cdots )$ &   & $7$ & $\infty$ & $\infty$ &3 & $\infty$ &$\infty$ \\
        & $(\pm t^2,\frac{1}{4}t^4,0,0,t,0,\cdots )$ &  &  &  &  & & & \\ \hline

     $(T_8)^{8}$ & $(t^3,0,-t^2,0,0,0,\cdots )$ &   & $\infty$ & $\infty$ & $\infty$ &$\infty$ & $\infty$ & $\infty$ \\
        & $(\pm t^2,0,0,0,t,0,\cdots )$ &  &  &  & & & &  \\ \hline

\end{tabular}

\smallskip

\caption{\small Lagrangian invariants for symplectic classes of $T_8$ singularity when $\omega\vert_ W\!=\! 0$; $W$ - the tangent space to a non-singular $3$-dimensional manifold  in $(\mathbb{R}^{2n\ge 6},\omega)$ containing $N\!\in \!(T_8)$.}\label{tabt8-lagr2}

\end{footnotesize}
\end{table}
\end{center}

\begin{rem}
We can see from Table \ref{tabt8-lagr2} that  the Lagrangian tangency order -- $Lt$ distinguishes different classes that the index of isotropy -- $ind$. For example the class $(T_8)^4$ in the case $c_1=0, c_2\ne 0$ and the class $(T_8)^{6,2}$ are distinguished by the index of isotropy - $ind$ but are not distinguished by the Lagrangian tangency order. We can distinguish these classes using the relative Lagrangian tangency order: for the class $(T_8)^4$ in the case $c_1=0, c_2\ne 0$ we have $Lt[C_2:C_1]=3$ and for the class $(T_8)^{6,2}$ we have $Lt[C_2:C_1]=4$.

  The index of isotropy -- $ind$ for the classes $(T_8)^{3,1}$, $(T_8)^4$, $(T_8)^{6,1}$, $(T_8)^{5,1}$ is less than for the class $(T_8)^{6,2}$ but the analogical inequality is not hold for the Lagrangian tangency order -- $Lt$.

We are not able to distinguish all symplectic classes using the Lagrangian tangency orders or the indices of isotropy but we can do it checking geometric conditions formulated analogically as for $T_7$ singularity (Section \ref{t7-geom_cond}).

\end{rem}

\section{Symplectic $T_7$-singularities}\label{sec-t7}

 Denote by $(T_7)$  the class of varieties in  a fixed symplectic space $(\mathbb R^{2n}, \omega )$ which are diffeomorphic to
\begin{equation}
\label{deft7} T_7=\{x\in \mathbb R ^{2n\geq
4}\,:x_1^2+x_2^3+x_3^3=x_2 x_3=x_{\geq 4}=0\}.\end{equation}

This is the classical $1$-dimensional isolated complete intersection singularity $T_7$ (\cite{G}, \cite{AVG}).  $N$ is quasi-homogeneous with weights $w(x_1)=3,\; w(x_2)=w(x_3)=2$.

We use the method of algebraic restrictions to obtain a complete classification of symplectic singularities of $(T_7)$ presented in the following theorem.

\begin{thm}\label{T7-main}
Any stratified submanifold of the symplectic space $(\mathbb R^{2n},\sum_{i=1}^n dp_i \wedge dq_i)$ which is diffeomorphic to $T_7$ is symplectically equivalent to one and only one of the normal forms $T_7^i, i = 0,1,\cdots ,7$ (resp. $i= 0,1,2,4$). The parameters $c, c_1, c_2$ of the normal forms are moduli.

\noindent $T_7^0$: \ $p_1^2 + p_2^3 + q_2^3 = 0, \ \ p_2q_2 = 0, \
\ q_1 = c_1q_2 + c_2p_2, \ \ p_{\ge 3} = q_{\ge 3} = 0, \ \ c_1\cdot c_2\ne 0$;

\smallskip

\noindent $T_7^1$: \ $p_1^2 + p_2^3 + q_1^3 = 0, \ \ p_2q_1 = 0, \
\ q_2 = c_1q_1 - c_2p_1p_2, \ \ p_{\ge 3} = q_{\ge 3} = 0$;

\smallskip

\noindent $T_7^2$: \ $p_1^2 + p_2^3 + q_2^3 = 0, \ \ p_2q_2 = 0, \
\ q_1 = \frac{c_1}{2}q_2^2 + \frac{c_2}{2}p_2^2, \ \ p_{\ge 3} =
q_{\ge 3} = 0, \ \ (c_1,c_2)\ne (0,0)$;

\smallskip

\noindent $T_7^4$: \ $p_1^2 + p_2^3 + q_2^3 = 0, \ \ p_2q_2 = 0, \
\ q_1 = \frac{c}{3}q_2^3, \ \ p_{\ge 3} = q_{\ge 3} = 0$;

\smallskip

\noindent $T_7^3$: \ $p_1^2 + p_2^3 + p_3^3 = 0, \  p_2p_3 = 0, \
q_1 = \frac{c_1}{2}p_2^2+\frac{1}{2}p_3^2,  \ q_2=-c_2p_1 p_3, \ p_{\ge 4} = q_{\ge 3} = 0;$

\smallskip

\noindent $T_7^5$:  \ $p_1^2 + p_2^3 + p_3^3 = 0, \  p_2p_3 = 0, \
q_1 = \frac{c}{3}p_3^3,  \ q_2=-p_1 p_3, \ p_{\ge 4} = q_{\ge 3} = 0;$

\smallskip

\noindent $T_7^6$:  \ $p_1^2 + p_2^3 + p_3^3 = 0, \  p_2p_3 = 0, \
q_1 = \frac{1}{3}p_3^3,   \ p_{\ge 4} = q_{\ge 2} = 0;$
\smallskip

\noindent $T_7^7$: \ $p_1^2 + p_2^3 + p_3^3 = 0, \ p_2p_3 = 0, \ q_{\ge 1} = p_{\ge 4} = 0.$
\end{thm}

 In Section \ref{t7-lagr} we use the Lagrangian tangency orders to distinguish more symplectic
singularity classes. In Section \ref{t7-geom_cond} we propose a geometric description of these singularities which confirms this more detailed classification. Some of the proofs are presented in
Section \ref{proofs}.

\subsection{Distinguishing symplectic classes of $T_7$ by Lagrangian tangency orders and the indices of isotropy}
\label{t7-lagr}

A curve $N\in (T_7)$ can be described as a union of two parametrical branches $B_1$ and $B_2$. Their
parameterization is given in the second column of Table \ref{tabt7-lagr-ind}. To distinguish the classes of this singularity  we need the following three invariants:
\begin{itemize}
  \item $Lt(N)=Lt(B_1,B_2)=\max\limits _L (\min \{t(B_1,L),t(B_2,L)\})$
  \item $L_{n}=\max\{Lt(B_1),Lt(B_2)\}=\max\{\max\limits _L\, t(B_1,L),\max\limits _L\, t(B_2,L)\}$
  \item $L_{f}=\min\{Lt(B_1),Lt(B_2)\}=\min\{\max\limits _L\, t(B_1,L),\max\limits _L\, t(B_2,L)\}$
\end{itemize}

where $L$ is a smooth Lagrangian submanifold of the symplectic space.

\medskip

Branches $B_1$ and $B_2$ are diffeomorphic and are not preserved by all symmetries of $T_7$ so  neither  $Lt(B_1)$ nor $Lt(B_2)$ can be used as invariants. The new invariants are defined instead:
$L_{n}$ describing the Lagrangian tangency order of the \emph{ nearest} branch and $L_{f}$ representing the Lagrangian tangency order of the \emph{ farthest } branch. Considering the triples
$(Lt(N), L_{n}, L_{f})$ we obtain more detailed classification of symplectic singularities of $T_7$ than the classification given in Table \ref{tabt7}. Some subclasses appear having a natural geometric interpretation (Tables \ref{tabt7-geom1} and \ref{tabt7-geom2}).

\begin{rem}
We can define the indices of isotropy for branches similarly as the Lagrangian tangency orders and use them to characterize singularities of $T_7$. We use the following invariants:
\begin{itemize}
    \item $ind_{n}=\max\{ind(B_1),ind(B_2)\}$
  \item $inf_{f}=\min\{ind(B_1),ind(B_2)\}$
\end{itemize}
where $ind(B_1), \ ind(B_2)$ denote the indices of isotropy for individual branches. They can be  calculated knowing their dependence on the Lagrangian tangency orders $Lt(B_1), Lt(B_2)$ for $A_2$ singularity (Table \ref{tab-a}).
\end{rem}

\begin{thm}
\label{lagr-t7} A stratified submanifold $N\in (T_7)$ of a symplectic space $(\mathbb R^{2n}, \omega )$ with the canonical coordinates $(p_1, q_1, \cdots, p_n, q_n)$ is symplectically equivalent to one and only one of the curves presented in the second column of Table \ref{tabt7-lagr-ind}. The parameters $c, c_1, c_2$ are moduli. The indices of isotropy are presented in the fourth, the fifth and the sixth column of Table \ref{tabt7-lagr-ind} and the Lagrangian tangency orders of the curve are presented in  the seventh, the eighth  and the ninth column of the Table.
\end{thm}

\begin{center}
\begin{table}[h]

    \begin{footnotesize}
    \noindent
    \begin{tabular}{|p{0.85cm}|p{3.4cm}|p{2.0cm}|c|c|c|c|c|c|}
                     \hline
    Class &  Parameterization \newline of branches  & Conditions \newline for subclasses  & $ind$ & $ind_n$ & $ind_f$ & $Lt(N)$   &  $L_{n}$ & $L_{f}$ \\ \hline
  $(T_7)^0$ \newline $2n\!\ge 4$   &  $(t^3,-c_1t^2,0,-t^2,0,\cdots )$
  \newline $(t^3,-c_2t^2,-t^2,0,0,\cdots )$ &  $c_1\cdot c_2\ne 0$ & $0$ & $0$ & $0$ & $2$ & $3$& $3$   \\  \hline

   &  &   $c_1\cdot c_2\ne 0$ & $0$ & $1$  & $0$ & $2$ & $5$& $3$ \\ \cline{3-9}
     $(T_7)^1$ & $(t^3,-t^2,0,-c_1t^2,0,\cdots )$ & $c_1=0, c_2\ne 0$& $0$ &  $1$ & $0$ & $3$ & $5$ & $3$ \\ \cline{3-9}
     $2n\!\ge 4$ & $(t^3,0,-t^2,c_2t^5,0,\cdots )$ & $c_1\ne 0, c_2=0$ & $0$ & $\infty$ & $0$ & $2$ & $\infty$ & $3$ \\ \cline{3-9}
           &   & $c_1=0, c_2=0$ & $0$ & $\infty$ & $0$ & $3$ & $\infty$ & $3$ \\ \hline

$(T_7)^2$  & $(t^3, \frac{c_1^2}{2}t^4, 0,-t^2, 0,\cdots)$ &  $c_1\cdot c_2\ne 0$ & $0$ & $1$ & $1$ & $2$ & $5$ & $5$   \\ \cline{3-9}
$2n\!\ge 4$  & $(t^3, \frac{c_2^2}{2}t^4,-t^2, 0, 0,\cdots)$ & $c_1\cdot c_2=\!0,\newline (c_1,c_2)\ne\!(0,0)$ & $0$ & $\infty$ & $1$ & $2$ & $\infty$ & $5$ \\ \hline

    $(T_7)^3$ & $(t^3\!, \frac{1}{2}t^4, 0, c_2t^5,-t^2, 0,\cdots\!)$ & $c_1\ne 0$ & $1$ & $1$ & $1$ & $5$ & $5$ & $5$ \\ \cline{3-9}
    $2n\!\ge 6$ & $(t^3, \frac{c_1}{2}t^4, -t^2, 0, 0, 0, \cdots)$ & $c_1=0$  & $1$ & $\infty$ & $1$ & $5$ & $\infty$ & $5$ \\ \hline

    $(T_7)^4$ \newline $2n\!\ge 4$ & $(t^3,\frac{c}{3}t^6,0,-t^2,0,\cdots )$
                            \newline $(t^3,0,-t^2,0,0,\cdots )$ & & $0$ & $\infty$ & $\infty$ & $2$ & $\infty$ & $\infty$  \\ \hline

    $(T_7)^5$ \newline $2n\!\ge 6$ & $(t^3\!,-\frac{c}{3}t^6, 0, t^5, -t^2, 0,\cdots\!)$
                        \newline $(t^3, 0, -t^2, 0, 0, 0,\cdots)$  & & $1$ & $\infty$ & $\infty$ & $5$ & $\infty$ & $\infty$ \\ \hline

     $(T_7)^6$ \newline $2n\!\ge 6$ & $(t^3\!,-\frac{1}{3}t^6, 0, 0, -t^2, 0,\cdots)$
                        \newline $(t^3, 0, -t^2, 0, 0, 0,\cdots)$  & & $2$ & $\infty$ & $\infty$ & $7$ & $\infty$ & $\infty$ \\ \hline
    $(T_7)^7$
    \newline $2n\!\ge 6$ & $(t^3, 0, 0, 0, -t^2, 0,\cdots)$
                 \newline $(t^3, 0, -t^2, 0, 0, 0,\cdots)$  & & $\infty$ & $\infty$ & $\infty$ & $\infty$ & $\infty$ & $\infty$ \\ \hline
\end{tabular}

\smallskip

\caption{\small  The Lagrangian tangency orders and the indices of isotropy for symplectic classes of $T_7$ singularity.}\label{tabt7-lagr-ind}

\end{footnotesize}
\end{table}
\end{center}
\medskip

The comparison of invariants presented in Table \ref{tabt7-lagr-ind} shows  that the Lagrangian tangency orders distinguish more symplectic classes than the indices of isotropy. The method of calculating these invariants is described in Section \ref{calc-lagr}.


\newpage

\subsection{Geometric conditions for the classes $(T_7)^i$}
\label{t7-geom_cond}

The classes $(T_7)^i$ can be distinguished geometrically, without using any local coordinate system.

\smallskip

Let $N\in (T_7)$. Then $N$ is the union of two branches -- singular $1$-dimensional irreducible components diffeomorphic to $A_2$ singularity. In local coordinates they have the form
\[\mathcal{B}_1=\{ x_1^2+x_3^3=0,\; x_2=x_{\geq 4}=0\},\]
\[\mathcal{B}_2=\{ x_1^2+x_2^3=0,\; x_{\geq 3}=0\}.\;\;\;\;\;\;\;\;\]
 Denote by $\ell _1, \ell _2$ the tangent lines at $0$ to the branches $\mathcal {B}_1$ and $\mathcal {B}_2$ respectively.
 These lines span a $2$-space $P_1$. Let $P_2$ be $2$-space tangent at $0$ to the branch $\mathcal {B}_1$ and $P_3$ be $2$-space tangent
 at $0$ to the branch $\mathcal {B}_2$. Define the line $\ell_3 = P_2\cap P_3$.  The lines $\ell _1, \ell _2, \ell _3$ span a  $3$-space $W=W(N)$.
Equivalently $W$ is  the tangent space at $0$ to some (and then any) non-singular $3$-manifold containing $N$.

\smallskip

The classes $(T_7)^i$ satisfy special conditions in terms of the restriction $\omega\vert_ W $, where $\omega $ is the symplectic form.
For $N=T_7=$(\ref{deft7}) it is easy to calculate
\begin{equation}
\label{lines} \ell _1 = \Span (\partial /\partial x_3),
\ \ell _2 = \Span (\partial /\partial x_2), \ \ell _3 = \Span (\partial
/\partial x_1).\end{equation}

\smallskip

\subsubsection{Geometric conditions for the class $[0]_{T_7}$}
\label{t7-zero_cond}

The geometric distinguishing of the class $(T_7)^7$ follows from Theorem \ref{thm B} :
$N\in (T_7)^7$ if and only if $N$ it is contained in a non-singular Lagrangian submanifold.
The following theorem gives a simple way to check the latter condition without using algebraic
restrictions. Given a $2$-form $\sigma $ on a non-singular submanifold $M$ of $\mathbb R^{2n}$
such that $\sigma (0)=0$ and a vector $v\in T_0M$ we denote by $\mathcal L_v\sigma $
the value at $0$ of the Lie derivative of $\sigma $ along a vector field $V$ on
$M$  such that $v = V(0)$.  The assumption $\sigma (0)=0$ implies that the choice of $V$ is irrelevant.

\begin{prop}
\label{check-zero-alg-rest} Let $N\in (T_7)$ be a stratified
submanifold of a symplectic space $(\mathbb R^{2n}, \omega )$. Let
$M^3$ be any non-singular submanifold containing $N$ and let
$\sigma $ be the restriction of $\omega $ to $TM^3$. Let $v_i\in
\ell _i$ be non-zero vectors. If the symplectic form $\omega $ has
zero algebraic restriction to $N$ then the following conditions are satisfied:

\smallskip

\noindent I.  \;$\sigma (0)=0$,
\smallskip

\noindent II. \;$\mathcal L_{v_3}\sigma (v_i, v_j) =0$ for $i,j\in\{1,2\}$,

\smallskip

\noindent III. \;$ \mathcal L_{v_i}\sigma (v_3, v_i) =0$ for $i\in\{1,2\}$,

\smallskip

\noindent IV.  \;$\mathcal L_{v_i}\sigma (v_3, v_j) = \mathcal L_{v_j}\sigma (v_3, v_i)$ for $i\neq j\in\{1,2\}$,

\end{prop}









\begin{thm}
\label{geom-cond-t7} A stratified submanifold $N\in (T_7)$ of a
symplectic space $(\mathbb R^{2n}, \omega )$ belongs to the class
$(T_7)^i$ if and only if the couple $(N, \omega )$ satisfies corresponding
conditions in the last column of Table \ref{tabt7-geom1} or \ref{tabt7-geom2}.
\end{thm}
\begin{center}
\begin{table}[h]
    \begin{small}
    \noindent
    \begin{tabular}{|p{0.7cm}|p{3.7cm}|p{7cm}|}
            \hline
    Class &  Normal form & Geometric conditions  \\ \hline

  $(T_7)^0$   &  $[T_7]^0: [\theta _1 + c_1\theta _2 + c_2\theta _3]_{T_7}$
                                        \newline  $c_1\cdot c_2\ne 0$  &  $  \omega|_{ \ell_i+\ell_j} \neq 0\;\; \forall i,j\in \{1,2,3\}$ so \newline 2-spaces tangent to branches are not isotropic   \\  \hline 

 $(T_7)^1$ & & $\exists i\!\ne\! j\in\! \{1,2\}\;\; \omega|_{\ell_i+\ell_3} =0$ and $\omega |_{\ell_j+\ell_3}\ne 0$ \newline
                    (exactly one branch has tangent 2-space  isotropic) \\ \cline{2-3}
         & $[T_7]^1_{2,5}:[c_1\theta _1 + \theta _2 + c_2\theta _5]_{T_7}$
   \newline $c_1\cdot c_2\ne 0$ & $\omega|_{\ell_1+\ell_2} \ne 0$ and  no branch is contained in a Lagrangian submanifold \\ \cline{2-3}
     & $[T_7]^1_{3,5}: [\theta _2 + c_2\theta _5]_{T_7}$,\newline $ c_2\ne 0$ & $\omega |_{\ell_1+\ell_2}=0$ and no branch is contained in a Lagrangian submanifold \\ \cline{2-3}
     & $[T_7]^1_{2,\infty}: [c_1\theta _1 + \theta _2 ]_{T_7}$,\newline $ c_1\ne 0$ &
                $\omega|_{\ell_1+\ell_2}\ne 0$ and exactly one branch is contained in a Lagrangian submanifold  \\ \cline{2-3}
                & $[T_7]^1_{3,\infty}: [ \theta _2 ]_{T_7}$ & $\omega|_{\ell_1+\ell_2}=0$ and exactly one branch is contained in a Lagrangian submanifold  \\ \hline 

$(T_7)^2$ & & $\omega |_{\ell_1+\ell_2}\neq 0, \omega|_{\ell_i+\ell_3} = 0\; \forall i \in\{1,2\}$\\ \cline{2-3}
  & $[T_7]^2_5: [\theta _1 + c_1\theta_4+c_2\theta_5]_{T_7}$
                                                    \newline  $c_1\cdot c_2\ne 0$ &   no branch is contained in a Lagrangian submanifold  \\ \cline{2-3}
 & $[T_7]^2_{\infty}: [\theta _1 + c_1\theta_4+c_2\theta_5]_{T_7}$
                                                    \newline  $c_1\cdot c_2= 0, \; c_1+c_2\ne 0$ &   exactly one branch is contained in a Lagrangian submanifold   \\ \hline 

    $(T_7)^4$  & $[T_7]^4: [\theta _1 + c\theta _7]_{T_7}$ &  $\omega |_{\ell_1+\ell_2}\neq 0, \omega |_{\ell_i+\ell_3} = 0\; \forall i \in\{1,2\}$,  \newline and  branches are contained in different Lagrangian submanifolds   \\ \hline
    \end{tabular}

\smallskip

\caption{\small Geometric interpretation of singularity classes of $T_7$ when $\omega\vert_ W \ne 0$;
$W$ - the tangent space to a non-singular $3$-dimensional manifold  in $(\mathbb{R}^{2n\ge 4},\omega)$ containing $N\in (T_7)$.}\label{tabt7-geom1}

\end{small}
\end{table}
\end{center}

\begin{center}
\begin{table}[h]

    \begin{small}
    \noindent
    \begin{tabular}{|p{0.8cm}|p{3.5cm}|p{7cm}|}
            \hline
    Class &  Normal form & Geometric conditions  \\ \hline

     $(T_7)^3$ & $[T_7]^3_5: [\theta _4 + c_1\theta _5+c_2\theta_6]_{T_7}$
            \newline $c_1\ne 0$ & III is not satisfied and no branch is contained in a Lagrangian submanifold \\ \cline{2-3}
     & $[T_7]^3_{\infty}: [\theta _4 +c_2\theta_6]_{T_7}$  &  III is not satisfied and exactly one branch is contained in a Lagrangian submanifold  \\ \hline

    $(T_7)^5$  & $[T_7]^5: [\theta _6 + c\theta _7]_{T_7}$ &  III is satisfied but  II is not
                                        and branches are contained in different  Lagrangian submanifolds. \\ \hline

     $(T_7)^6$  & $[T_7]^6: [\theta _7]_{T_7}$  & I - IV are satisfied and
                                          branches are contained in different  Lagrangian submanifolds.   \\ \hline
    $(T_7)^7$ & $[T_7]^7: [0]_{T_7}$ &    I - IV are satisfied and
            $N$  is contained in a  Lagrangian submanifold   \\ \hline
\end{tabular}

\smallskip

\caption{\small Geometric interpretation of singularity classes of $T_7$ when $\omega\vert_ W\!=\! 0$; $W$ - the tangent space to a non-singular $3$-dimensional manifold  in $(\mathbb{R}^{2n\ge 6},\omega)$ containing $N\!\in \!(T_7)$;
 I\! -\! IV\! -- conditions of Proposition \ref{check-zero-alg-rest}.}\label{tabt7-geom2}

\end{small}
\end{table}
\end{center}
The proofs of the theorems of this Section are presented in Section \ref{t7-geom}.

\section{Proofs}\label{proofs}

\subsection{The method of algebraic restrictions}
\label{method}

In this section we present basic facts on the method of algebraic restrictions, which is a very powerful tool for the symplectic classification. The details of the method and proofs of all results of this section can be found in \cite{DJZ2}.

  Given a germ of a non-singular manifold $M$ denote by $\Lambda ^p(M)$ the space of all germs at $0$ of differential $p$-forms on $M$. Given a subset $N\subset M$ introduce the following subspaces of $\Lambda ^p(M)$:
$$\Lambda ^p_N(M) = \{\omega \in \Lambda ^p(M): \ \ \omega (x)=0 \ \text {for any} \ x\in N \};$$
$$\mathcal A^p_0(N, M) = \{\alpha  + d\beta : \ \ \alpha \in \Lambda _N^p(M), \ \beta \in \Lambda _N^{p-1}(M).\}$$

\smallskip

\begin{defn}
\label{main-def} Let $N$ be the germ of a subset of $M$ and let
$\omega \in \Lambda ^p(M)$. The {\bf algebraic restriction} of
$\omega $ to $N$ is the equivalence class of $\omega $ in $\Lambda
^p(M)$, where the equivalence is as follows: $\omega $ is
equivalent to $\widetilde \omega $ if $\omega - \widetilde \omega
\in \mathcal A^p_0(N, M)$.
\end{defn}

\noindent {\bf Notation}. The algebraic restriction of the germ of
a $p$-form $\omega $ on $M$ to the germ of a subset $N\subset M$
will be denoted by $[\omega ]_N$. Writing $[\omega ]_N=0$ (or
saying that $\omega $ has zero algebraic restriction to $N$) we
mean that $[\omega ]_N = [0]_N$, i.e. $\omega \in A^p_0(N, M)$.

\medskip




\begin{defn}Two algebraic restrictions
$[\omega ]_N$ and $[\widetilde \omega ]_{\widetilde N}$ are called {\bf
diffeomorphic} if there exists the germ of a diffeomorphism $\Phi:
\widetilde M\to M$ such that $\Phi(\widetilde N)=N$ and  $\Phi ^*([\omega ]_N) =[\widetilde \omega ]_{\widetilde N}$.
\end{defn}

\smallskip


The method of algebraic restrictions applied to singular
quasi-homogeneous subsets is based on the following theorem.

\begin{thm}[Theorem A in \cite{DJZ2}] \label{thm A}
Let $N$ be the germ of a quasi-homogeneous subset of $\mathbb
R^{2n}$. Let $\omega _0, \omega _1$ be germs of symplectic forms
on $\mathbb R^{2n}$ with the same algebraic restriction to $N$.
There exists a local diffeomorphism $\Phi $ such that $\Phi (x) =
x$ for any $x\in N$ and $\Phi ^*\omega _1 = \omega _0$.

Two germs of quasi-homogeneous subsets $N_1, N_2$ of a fixed
symplectic space $(\mathbb R^{2n}, \omega )$ are symplectically
equivalent if and only if the algebraic restrictions of the
symplectic form $\omega $ to $N_1$ and $N_2$ are diffeomorphic.

\end{thm}

\medskip

Theorem \ref{thm A} reduces the problem of symplectic
classification of germs of singular quasi-homogeneous subsets to
the problem of diffeomorphic classification of algebraic
restrictions of the germ of the symplectic form to the germs of
singular quasi-homogeneous subsets.

The geometric meaning of zero algebraic restriction is explained
by the following theorem.

\begin{thm}[Theorem {\bf B} in \cite{DJZ2}] \label{thm B}  {\it The germ of a quasi-homogeneous
 set $N$  of a symplectic space
$(\mathbb R^{2n}, \omega )$ is contained in a non-singular
Lagrangian submanifold if and only if the symplectic form $\omega
$ has zero algebraic restriction to $N$.}\
\end{thm}


The following result shows that the method of algebraic
restrictions is very powerful tool in symplectic classification of
singular curves.

\begin{thm}[Theorem 2 in \cite{D}]
\label{main-alg} Let $C$ be the germ of a $\mathbb K$-analytic
curve (for $\mathbb K=\mathbb R$ or $\mathbb K=\mathbb C$). Then
the space of algebraic restrictions of germs of closed $2$-forms
to $C$ is a finite dimensional vector space.
\end{thm}

By a {\bf $\mathbb K$-analytic curve} we understand a subset of
$\mathbb K^m$ which is locally diffeomorphic to a $1$-dimensional
(possibly singular) $\mathbb K$-analytic subvariety of $\mathbb
K^m$. Germs of $\mathbb C$-analytic parameterized curves can be
identified with germs of irreducible $\mathbb C$-analytic curves.

In the paper we use the following notations:

\smallskip

\noindent $\bullet$  $\algrestall $: \ \ the vector space
consisting of algebraic restrictions of germs of all $2$-forms on
$\mathbb R^{2n}$ to the germ of a subset $N\subset \mathbb
R^{2n}$;

\smallskip

\noindent $\bullet$  $\algrestclosed$: \ \ the subspace of
$\algrestall $ consisting of algebraic restrictions of germs of
all closed $2$-forms on $\mathbb R^{2n}$ to $N$;

\smallskip

\noindent $\bullet$ $\algrest $: \ \ the open set in
$\algrestclosed$ consisting of algebraic restrictions of germs of
all symplectic $2$-forms on $\mathbb R^{2n}$ to $N$.

\medskip

For calculating discrete invariants we use the following propositions.

\begin{prop}[\cite{DJZ2}]\label{sm}
The symplectic multiplicity  of the germ of a quasi-homogeneous subset $N$ in a symplectic space is equal to the codimension of the orbit of the algebraic restriction $[\omega ]_N$ with respect to the group of local diffeomorphisms preserving $N$  in the space of algebraic restrictions of closed  $2$-forms to $N$.
\end{prop}

\begin{prop}[\cite{DJZ2}]\label{ii}
The index of isotropy  of the germ of a quasi-homogeneous subset $N$ in a symplectic space $(\mathbb R^{2n}, \omega )$ is equal to the maximal order of vanishing of closed $2$-forms representing the algebraic restriction $[\omega ]_N$.
\end{prop}

\begin{prop}[\cite{D}]\label{lto}
Let $f$ be the germ of a quasi-homogeneous curve such that the algebraic restriction of a symplectic form to it can be represented by a closed $2$-form vanishing at $0$. Then the Lagrangian tangency order of the germ of a quasi-homogeneous curve $f$ is the maximum of the order of vanishing on $f$ over all $1$-forms $\alpha$ such that $[\omega]_f=[d\alpha]_f$
\end{prop}

\subsection{Algebraic restrictions to $T_7$ and their classification}\label{t7-class}

One has the following relations for $(T_7)$-singularities
\begin{equation}
[d(x_2 x_3)]_{T_7}=[x_2 dx_3+x_3 dx_2]_{T_7}=0
\label{t1}
\end{equation}
\begin{equation}
[d(x_1^2+x_2^3+x_3^3)]_{T_7}=[2x_1dx_1+3x_2^2dx_2+3x_3^2dx_3]_{T_7}=0
\label{t2}
\end{equation}
Multiplying these relations by suitable $1$-forms we obtain the relations in Table \ref{tabTk1}.

\renewcommand*{\arraystretch}{1.3}
\begin{small}
\begin{table}[h]
\begin{center}
\begin{tabular}{|c|c|c|}

 \hline

        & relations & proof\\ \hline

   1. & $[x_2dx_2 \wedge dx_3]_N=0$ & (\ref{t1})$\wedge\, dx_2$ \\ \hline

   2. &  $[x_3dx_2 \wedge dx_3]_N=0$ & (\ref{t1})$\wedge\, dx_3$\\ \hline

    3. &  $[x_3dx_1 \wedge dx_2]_N=[x_2dx_3\wedge dx_1]_N$ & (\ref{t1})$\wedge\, dx_1$\\\hline

    4. &  $[x_1dx_1 \wedge dx_2]_N=0$  & (\ref{t2})$\wedge\, dx_2\;$ and row 2. \\\hline

    5. & $[x_1dx_1 \wedge dx_3]_N=0$  & (\ref{t2})$\wedge\, dx_3\;$ and row 1. \\\hline

    6. &  $[x_2^2dx_1 \wedge dx_2]_N=[x_3^{2}dx_3\wedge dx_1]_N$ & (\ref{t2})$\wedge\, dx_1$ \\ \hline

    7. & $[x_1^2dx_2 \wedge dx_3]_N=0$ &  rows 1. and 2.   and $[x_1^2]_N=[-x_2^3-x_3^3]_N$ \\ \hline

    8. & $[x_3^2dx_1 \wedge dx_2]_N=0$ & (\ref{t1})$\wedge\, x_3 dx_1$  and $[x_2 x_3]_N=0$\\  \hline

\end{tabular}
\end{center}
\smallskip
\caption{\small Relations towards calculating $[\Lambda^2(\mathbb R^{2n})]_N$ for $N=T_7$}\label{tabTk1}
\end{table}
\end{small}
\smallskip

 Using the method of algebraic restrictions and Table \ref{tabTk1} we obtain the following proposition:

\begin{prop}
\label{Tk-all}
$[\Lambda ^{2}(\mathbb R^{2n})]_{T_7}$ is a $8$-dimensional vector space spanned by the algebraic restrictions to $T_7$ of the $2$-forms

\smallskip

$\theta _1= dx_2\wedge dx_3, \;\; \theta _2=dx_1\wedge dx_3,\;\; \theta_3 = dx_1\wedge dx_2,$

\smallskip

  $\theta _4 = x_3dx_1\wedge dx_3,\;\; \theta _5 = x_2dx_1\wedge dx_2,$

\smallskip

$\sigma _1 = x_3dx_1\wedge dx_2,\;\;  \sigma _2 = x_1 dx_2\wedge dx_3, $

\smallskip

$\theta _7= x_3^2 dx_1\wedge dx_3$.

\end{prop}

Proposition \ref{Tk-all} and results of Section \ref{method}  imply the following description of the space $[Z ^2(\mathbb R^{2n})]_{T_7}$ and the manifold $[{\rm Symp} (\mathbb R^{2n})]_{T_7}$.

\begin{thm} \label{t7-baza}
 $[Z^2(\mathbb R^{2n})]_{T_7}$ is a $7$-dimensional vector space
 spanned by the algebraic restrictions to $T_7$  of the quasi-homogeneous $2$-forms $\theta_i$ 

\smallskip

$$\theta_1,\,\theta_2,\,\theta_3,\,\theta_4,\,\theta_5,\,\theta_6=\sigma_1-\sigma_2,\,\theta_7.$$

If $n\ge 3$ then $[{\rm Symp} (\mathbb R^{2n})]_{T_7} = [Z^2(\mathbb R^{2n})]_{T_7}$. The manifold $[{\rm Symp} (\mathbb R^{4})]_{T_7}$ is an open part of the $7$-space $[Z^2 (\mathbb R^{4})]_{T_7}$ consisting of algebraic restrictions of the form $[c_1\theta _1 + \cdots + c_7\theta _7]_{T_7}$ such that $(c_1,c_2,c_3)\ne (0,0,0)$.
\end{thm}

\begin{thm}
\label{klasT7} $ \ $

\smallskip

\noindent (i) \ Any algebraic restriction in $[Z ^2(\mathbb R^{2n})]_{T_7}$ can be brought by a symmetry of $T_7$ to one of the normal forms $[T_7]^i$ given in the second column of Table \ref{tabt7};

\smallskip

\noindent (ii)  \ The codimension in $[Z ^2(\mathbb R^{2n})]_{T_7}$ of the singularity class corresponding to the normal form $[T_7]^i$ is equal to $i$;

\smallskip

\noindent (iii) \ The singularity classes corresponding to the normal forms are disjoint;

\smallskip

\noindent (iv) \ The parameters $c, c_1, c_2$ of the normal forms $[T_7]^0, [T_7]^1, [T_7]^2, [T_7]^3, [T_7]^4, [T_7]^5$ are moduli.

\end{thm}

\renewcommand*{\arraystretch}{1.3}
\begin{center}
\begin{table}[h]

    \begin{small}
    \noindent
    \begin{tabular}{|p{2.5cm}|p{6.0cm}|c|c|c|}
                 \hline
    Symplectic class &   Normal forms for algebraic restrictions    & $cod$ & $\mu ^{\rm sym}$ &  $ind$  \\ \hline
  $(T_7)^0$ \;\;  $(2n\ge 4)$ & $[T_7]^0: [\theta _1 + c_1\theta _2 + c_2\theta _3]_{T_7}$,\;
                         $c_1\cdot c_2\ne 0$ &  $0$ & $2$ & $0$  \\  \hline
  $(T_7)^1$ \;\; $(2n\ge 4)$ & $[T_7]^1: [c_1\theta _1 + \theta _2 + c_2\theta _5]_{T_7}$
                                        &  $1$ & $3$ & $0$ \\ \hline
    $(T_7)^2$ \;\; $(2n\ge 4)$& $[T_7]^2: [\theta _1 + c_1\theta_4+c_2\theta_5]_{T_7}$,\;
            $(c_1,c_2)\ne (0,0)$  & $2$ & $4$ & $0$     \\ \hline
    $(T_7)^3$ \;\; $(2n\ge 6)$ & $[T_7]^3: [\theta _4 + c_1\theta _5+c_2\theta_6]_{T_7}$ &
                $3$ & $5$ & $1$  \\ \hline
    $(T_7)^4$ \;\; $(2n\ge 4)$ & $[T_7]^4: [\theta _1 + c\theta _7]_{T_7}$
                                         &  $4$ & $5$ & $0$    \\ \hline
    $(T_7)^5$ \;\; $(2n\ge 6)$ & $[T_7]^5: [\theta _6 + c\theta _7]_{T_7}$
                                         &  $5$ & $6$ & $1$ \\ \hline
     $(T_7)^6$ \;\; $(2n\ge 6)$ & $[T_7]^6: [\theta _7]_{T_7}$
                                         &  $6$ & $6$ & $2$   \\ \hline
    $(T_7)^7$ \;\; $(2n\ge 6)$ & $[T_7]^7: [0]_{T_7}$ &  $7$ & $7$ & $\infty $ \\ \hline
\end{tabular}

\smallskip

\caption{\small Classification of symplectic $T_7$ singularities.  \newline
$cod$ -- codimension of the classes; \ $\mu ^{sym}$-- symplectic multiplicity; \newline $ind$ -- the index of isotropy.}\label{tabt7}

\end{small}
\end{table}
\end{center}

\medskip

The proof of Theorem \ref{klasT7} is presented in section \ref{t7-proof}.

\noindent In the first column of Table \ref{tabt7}  by $(T_7)^i$ we denote a subclass of $(T_7)$ consisting of $N\in (T_7)$ such that the algebraic restriction $[\omega ]_N$ is diffeomorphic to some algebraic restriction of the normal form $[T_7]^i$. Theorem \ref{thm A}, Theorem \ref{klasT7} and Proposition \ref{t7-baza} imply the following statement which  explains why the given stratification of $(T_7)$ is natural.

\begin{thm}
\label{t7-char-classes}  Fix $i \in \{0,1,\cdots,7\}$. All stratified submanifolds $N\in (T_7)^i$ have the same (a) symplectic multiplicity and (b) index of isotropy given in Table \ref{tabt7}.
\end{thm}

\begin{proof} Part (a) follows from Proposition  \ref{sm} and Theorem \ref{klasT7} and
the fact that the codimension in $[Z ^2(\mathbb R^{2n})]_{T_7}$ of the orbit of an algebraic restriction $a\in [T_7]^i$ is equal to the sum of the number of moduli in the normal form $[T_7]^i$ and the codimension in $[Z ^2(\mathbb R^{2n})]_{T_7}$ of the class of algebraic restrictions defined by this normal form.

Part (b) follows from Theorem \ref{thm B} and Proposition \ref{ii}.
\end{proof}

\begin{prop}
\label{def-classes-t7} The classes $(T_7)^i$ are symplectic singularity classes, i.e. they are closed with respect to the action of the group of symplectomorphisms. The class $(T_7)$ is the disjoint union of the classes $(T_7)^i, i \in \{0,1,\cdots,7\}$. The classes $(T_7)^0, (T_7)^1, (T_7)^2, (T_7)^4$ are non-empty for any dimension $2n\ge 4$ of the symplectic space; the classes $(T_7)^3,(T_7)^5, (T_7)^6, (T_7)^7$ are empty if $n=2$ and not empty if $n\ge 3$.
\end{prop}

\subsection{Symplectic normal forms. Proof of Theorem \ref{T7-main}}
\label{t7-normal}

Let us transfer the normal forms $[T_7]^i$  to symplectic normal forms
. Fix a family $\omega ^i$ of symplectic forms on $\mathbb R^{2n}$ realizing the family $[T_7]^i$ of algebraic restrictions.  We can fix, for example

\smallskip

\noindent $\omega ^0 = \theta _1 + c_1\theta _2 + c_2\theta _3 + dx_1\wedge dx_4 + dx_5\wedge dx_6 + \cdots + dx_{2n-1}\wedge dx_{2n}, \ \ c_1\cdot c_2\ne 0;$

\smallskip

\noindent $\omega ^1 = c_1 \theta _1 + \theta _2 + c_2\theta _5 + dx_2\wedge dx_4 + dx_5\wedge dx_6 + \cdots + dx_{2n-1}\wedge dx_{2n}; $

\smallskip

\noindent $\omega ^2 = \theta_1+ c_1 \theta _4 + c_2 \theta _5 + dx_1\wedge dx_4 + dx_5\wedge dx_6 + \cdots + dx_{2n-1}\wedge dx_{2n},\ \ (c_1,c_2)\ne (0,0);$

\smallskip

\noindent $\omega ^3 = \theta _4 + c_1\theta _5+c_2\theta_6 + dx_1\wedge dx_4 + dx_2\wedge dx_5 + dx_3\wedge dx_6 + dx_7\wedge dx_8 + \cdots + dx_{2n-1}\wedge dx_{2n};$

\smallskip

\noindent $\omega ^4 = \theta _1 + c\theta _7 + dx_1\wedge dx_4 + dx_5\wedge dx_6 + \cdots + dx_{2n-1}\wedge dx_{2n};$

\smallskip

\noindent $\omega ^5 = \theta _6 + c\theta _7+ dx_1\wedge dx_4 + dx_2\wedge dx_5 + dx_3\wedge dx_6 + dx_7\wedge dx_8 + \cdots + dx_{2n-1}\wedge dx_{2n};$

\smallskip

\noindent $\omega ^6 = \theta _7+ dx_1\wedge dx_4 + dx_2\wedge dx_5 + dx_3\wedge dx_6 + dx_7\wedge dx_8 + \cdots + dx_{2n-1}\wedge dx_{2n};$

\smallskip

\noindent $\omega ^7 = dx_1\wedge dx_4 + dx_2\wedge dx_5 + dx_3\wedge dx_6 + dx_7\wedge dx_8 + \cdots + dx_{2n-1}\wedge dx_{2n}.$

\smallskip

 Let $\omega = \sum_{i=1}^m dp_i \wedge dq_i$, where $(p_1,q_1,\cdots,p_n,q_n)$ is the coordinate system on $\mathbb R^{2n}, n\ge 3$ (resp. $n=2$). Fix, for $i=0,1,\cdots ,7$ (resp. for $i = 0,1,2,4)$ a family $\Phi ^i$ of local diffeomorphisms which bring the family of symplectic forms $\omega ^i$ to the symplectic form $\omega $: $(\Phi ^i)^*\omega ^i = \omega $. Consider the families $T_7^i = (\Phi ^i)^{-1}(T_7)$. Any stratified submanifold of the symplectic space $(\mathbb R^{2n}, \omega )$ which is diffeomorphic to $T_7$ is symplectically equivalent to one and only one of the normal forms $T_7^i, i = 0,1,\cdots ,7$ (resp. $i= 0,1,2,4$) presented in Theorem \ref{T7-main}. By Theorem \ref{klasT7} we obtain that  parameters $c,c_1,c_2$ of the normal forms are moduli.

\subsection{ Proof of Theorem \ref{t7-lagr}}
\label{calc-lagr}
The numbers $ind(B_1)$ and $ind(B_2)$ are computed using Proposition \ref{ii} to branches $B_1$ and $B_2$. The space $[Z^2 (\mathbb R^{2n})]_{B_1}$ is spanned only by the algebraic restrictions to $B_1$ of the $2$-forms $\theta_2, \,\theta_4$. The space $[Z^2 (\mathbb R^{2n})]_{B_2}$ is spanned only by the algebraic restrictions to $B_2$ of the $2$-forms $\theta_3, \,\theta_5$. Branches are curves of type $A_2$ and from Table \ref{tab-a} we know the interaction between the index of isotropy and the Lagrangian tangency order. Knowing $ind(B_1)$ and $ind(B_2)$ we obtain $Lt(B_1)=3+ind(B_1)$ and $Lt(B_2)=3+ind(B_2)$.  Then $L_{f}$ is the minimum of these numbers and $L_{n}$ is the maximum of them. Next we calculate $Lt(N)$ by definition finding the nearest Lagrangian submanifold to the branches knowing that it can not be greater than $L_{f}$.



As an example we calculate the invariants for the class $(T_7)^1$.

\noindent We have $[\omega^1]_{B_1}=[c_1 \theta_1+\theta_2+c_2\theta_5]_{B_1}=[\theta_2]_{B_1}$ and thus $ind(B_1)=0$ and $Lt(B_1)= 3$. 

\noindent $[\omega^1]_{B_2}=[c_1 \theta_1+\theta_2+c_2\theta_5]_{B_2}=[c_2\theta_5]_{B_2}$ and thus $ind(B_2)=1$ and $Lt(B_2)= 5$ if $c_2\ne 0$ and $ind(B_2)=\infty$ and $Lt(B_2)= \infty$ if $c_2= 0$.

\noindent Finally  for the class $(T_7)^1$ we have $L_{n}=5$ if $c_2\ne 0$
and $L_{n}=\infty$ if $c_2= 0$ and $L_{f}=3$ so $Lt(N)\leq 3$.

\noindent For the smooth Lagrangian submanifolds $L$ defined by the conditions: \par\noindent
$\;p_1=0,\; q_2=0$ and $ p_i=0$ for $ i> 2$ we get $t[N,L]=3$ if $c_1=0$ thus $Lt(N)= 3$ in this case. But if $c_1\ne 0$ then $t[N,L]=2$ and it can not be greater for any other smooth Lagrangian submanifold so $Lt(N)= 2$ in this case.

\subsection{ Proof of Theorem \ref{geom-cond-t7}}
\label{t7-geom}

\begin{proof}[Proof of Proposition \ref{check-zero-alg-rest}]
Any $2$-form $\sigma$ which has zero algebraic restriction to $T_7$  can be expressed in the following form $\sigma=H_1 \alpha+H_2 \beta+dH_1 \wedge \gamma+dH_2 \wedge \delta$, where $H_1=x_1^2+x_2^3+x_3^3$, $H_2=x_2 x_3$ and $\alpha, \beta$ are $2$-forms on $TM^3$ and $\gamma=\gamma_1 dx_1+\gamma_2 dx_2+\gamma_3 dx_3$ and $\delta=\delta_1 dx_1+\delta_2 dx_2+\delta_3 dx_3$ are $1$-forms on $TM^3$. Since
\begin{equation}\label{H0}
H_1(0)=H_2(0)=0, \ \ dH_1|_0=dH_2|_0=0
\end{equation}
we obtain the following equality
 $$\mathcal L_v\sigma = d(V\ \rfloor \sigma)|_0+(V\ \rfloor d\sigma)|_0 = d(V\ \rfloor  \sigma)|_0.$$
(\ref{H0}) also implies that
$$d(V\ \rfloor \sigma)|_0=d(V\rfloor dH_1)|_0\wedge \gamma|_0+d(V\rfloor dH_2)|_0\wedge \delta|_0.$$

By simply calculation we get

$$\mathcal L_{v_1}\sigma= dx_2\wedge\delta |_0=\delta_3|_0\, dx_2\wedge dx_3-\delta_1|_0\, dx_1\wedge dx_2,$$

$$\mathcal L_{v_2}\sigma= dx_3\wedge\delta |_0=\delta_1|_0\, dx_3\wedge dx_1-\delta_2|_0\, dx_2\wedge dx_3,$$

$$\mathcal L_{v_3}\sigma= 2dx_1\wedge\gamma |_0=2\gamma_2|_0\, dx_1\wedge dx_2-2\gamma_3|_0\, dx_3\wedge dx_1.$$

Finally we obtain

$$ \mathcal L_{v_1}\sigma (v_3, v_1) =0, \ \   \mathcal L_{v_2}\sigma (v_3, v_2) =0, \  \ \mathcal L_{v_3}\sigma (v_1, v_2) =0,$$
$$\mathcal L_{v_1}\sigma (v_3, v_2) =-\delta_1|_0= \mathcal L_{v_2}\sigma (v_3, v_1).$$

\end{proof}

\begin{proof}[Proof of Theorem \ref{geom-cond-t7}]
The conditions on the pair $(\omega, N)$ in the last column of Table \ref{tabt7-geom1} and Table \ref{tabt7-geom2} are disjoint. It suffices to prove that these conditions  the row of $(T_7)^i$, are satisfied for any $N\in (T_7)^i$. This  is a corollary of the following claims:

\smallskip

\noindent 1. Each of the conditions in the last column of Tables \ref{tabt7-geom1}, \ref{tabt7-geom2}
is invariant with respect to the action of the group of diffeomorphisms in the space of pairs $(\omega , N)$;

\smallskip

\noindent 2. Each of these conditions depends only on the algebraic restriction $[\omega ]_N$;

\smallskip

\noindent 3. Take the simplest $2$-forms $\omega ^i$ representing
the normal forms $[T_7]^i$ for algebraic restrictions:
 $\omega ^0, \ \omega ^1, \ \omega ^2, \ \omega ^3, \ \omega ^4, \ \omega ^5, \ \omega ^6, \ \omega ^7.$
 The pair $(\omega = \omega ^i, T_7)$ satisfies the condition
in the last column of Table \ref{tabt7-geom1} or Table \ref{tabt7-geom2}, the row of $(T_7)^i$.

\medskip


To prove the third statement we note that in the case
$N = T_7 = (\ref{deft7})$ one has $W = \Span
(\partial/\partial x_1, \partial / \partial x_2, \partial / \partial x_3)$ and $v_1\in \ell_1= \Span (\partial / \partial x_3)$, $v_2\in \ell_2=span(\partial / \partial x_2)$, $v_3\in \ell_3= \Span(\partial / \partial x_1)$. By simply calculation and observation of Lagrangian tangency orders we obtain that following statements are true:
\smallskip

\noindent $(T^0)$ $\omega^0|_{\ell_1+\ell_2}\ne 0$ and
$\omega^0|_{\ell_1+\ell_3}\ne 0$ and also
$\omega^0|_{\ell_2+\ell_3}\ne 0$, and $L_{n}< \infty$ and $L_{f}<
\infty$ hence no branch is contained in a smooth Lagrangian
submanifold.

\smallskip

\noindent $(T^1)$ For any $c_1, c_2$\; $\omega^1|_{\ell_1+\ell_3}=
0$ and $\omega^1|_{\ell_2+\ell_3}\ne 0$ or
$\omega^1|_{\ell_1+\ell_3}\ne 0$ and
$\omega^1|_{\ell_2+\ell_3}=0$. If $c_2=0$ then  and $L_{n}=
\infty$ and $L_{f}< \infty$ hence exactly one  branch is contained
in some smooth Lagrangian submanifold. For $c_2\ne 0$ $L_{n}<
\infty$ and $L_{f}< \infty$ so no branch is contained in a smooth
Lagrangian submanifold. $\omega^1|_{\ell_1+\ell_2}= 0$ if and only
if $c_1=0$.

\smallskip

\noindent $(T^2)$ For any $c_1, c_2$\;
$\omega^2|_{\ell_1+\ell_2}\ne 0$ and $\omega^2|_{\ell_1+\ell_3}=
0$ and also $\omega^2|_{\ell_2+\ell_3}= 0$. If $c_1\cdot c_2\ne 0$
then $L_{n}< \infty$ and $L_{f} < \infty$ so no branch is
contained in a Lagrangian submanifold. If $c_1= 0$ and $c_2\ne 0$
or $c_1\ne 0$ and $c_2= 0$ then and $L_{n}= \infty$ and $L_{f}<
\infty$ hence exactly one  branch is contained in some smooth
Lagrangian submanifold.

\smallskip

\noindent $(T^3)$ The  Lie derivative of $\omega ^3= \theta _4 +
c_1\theta _5+c_2\theta_6 $ along a vector field $V=\partial
/\partial x_3 $ is not equal to $0$, so condition III of
Proposition \ref{check-zero-alg-rest} is not satisfied. If $c_1\ne
0$ then $L_{n}< \infty$ and $L_{f}< \infty$ hence no branch is
contained in a Lagrangian submanifold. If $c_1= 0$  then $L_{n}=
\infty$ and $L_{f}< \infty$ hence only one branch is contained in
some Lagrangian submanifold.

\smallskip

\noindent $(T^4)$ For any $c$\; $\omega^4|_{\ell_1+\ell_2}\ne 0$
and $\omega^4|_{\ell_1+\ell_3}= 0$ and also
$\omega^4|_{\ell_2+\ell_3}= 0$. Both branches are contained in
different Lagrangian submanifolds since $L_{n}= L_{f}=\infty$ and
$Lt(N)<\infty$.

\smallskip

\noindent $(T^5)$  We can calculate the  Lie derivatives of
$\omega ^5= \theta _6 + c\theta _7 $ along a vector fields
$V_1=\partial /\partial x_3 $ and $V_2=\partial /\partial x_2 $
and $V_3=\partial /\partial x_3 $: $ \mathcal L_{V_1}\omega ^5
(V_3, V_1) =0$ and $ \mathcal L_{V_2}\omega ^5 (V_3, V_2) =0$, so
 condition III of Proposition \ref{check-zero-alg-rest} is
satisfied, but the  Lie derivative
 $\mathcal L_{V_3}\omega ^5 (V_1, V_2)$ is not equal to $0$, so  condition II of  Proposition \ref{check-zero-alg-rest} is not satisfied.
 We have $Lt(N)< \infty$ and $L_{n}=L_{f}=\infty$ hence  branches are contained in different Lagrangian submanifolds.

 \smallskip

\noindent $(T^6)$ The Lie derivatives of $\omega ^6 = \theta _7$,
$ \mathcal L_{V_i}\omega ^6 (V_j, V_k) =0$ for $i,j,k\in
\{1,2,3\}$, so  conditions II, III and IV of Proposition
\ref{check-zero-alg-rest} are satisfied. We have $Lt(N)< \infty$
and $L_{n}=L_{f}=\infty$ hence  branches are contained in
different Lagrangian submanifolds.

\smallskip

\noindent $(T^7)$ For $\omega ^7 = 0$ we have $ \mathcal
L_{V_i}\omega ^7 (V_j, V_k) =0$ for $i,j,k\in \{1,2,3\}$, so
conditions II, III and IV of  Proposition
\ref{check-zero-alg-rest} are satisfied. The condition
$Lt(N)=\infty$ implies the curve N is contained in a smooth
Lagrangian submanifold.

\end{proof}

\subsection{Proof of Theorem \ref{klasT7}}
\label{t7-proof}

In our proof we use vector fields tangent to $N\in (T_7)$. A Hamiltonian vector field is an example of such a vector field. We recall by \cite{ArGorLyVa} a suitable definition and facts.

\begin{defn}
 Let $H=\left\{H_1=\cdots=H_p=0\right\}\subset \mathbb R^n$ be a
complete intersection. Consider a set of $p+1$ integers $1\leq
i_1<\cdots<i_{p+1}\leq n$. A Hamiltonian vector field
$X_H(i_1,\ldots,i_{p+1})$ on a complete intersection
 $H$ is the determinant obtained by expansion with respect to the first row of the symbolic $(p+1)\times(p+1)$ matrix

\begin{equation}  \label{hamiltonfield}
X_H(i_1,\ldots,i_{p+1})= det
\left[ \begin{array}{ccc}
\partial /\partial x_{i_1} & \cdots & \partial /\partial x_{i_{p+1}} \\
\partial H_1 /\partial x_{i_1} & \cdots & \partial H_1 /\partial x_{i_{p+1}} \\
\vdots &    \ldots   &  \vdots \\
\partial H_p /\partial x_{i_1} & \cdots & \partial H_p /\partial x_{i_{p+1}}
\end{array} \right]
\end{equation}

\end{defn}

\begin{thm} [\cite{Wa}]
Let $H=\left\{H_1=\cdots=H_p=0\right\}\subset \mathbb R^n$ be a
positive dimensional complete intersection with an isolated
singularity. If $H_1,\dots,H_p$ are quasi-homoge\-neous with
positive weights $\lambda_1,\ldots,\lambda_n$ than the module of
vector fields tangent to $H$  is generated by the Euler vector
field $E=\sum_{i=1}^n\lambda_i x_i \frac{\partial}{\partial x_i}$
and the Hamiltonian fields $X_H(i_1,\ldots,i_{p+1})$ where the
numbers $i_1,\ldots,i_{p+1}$ run through all possible sets $1\leq
i_1<\cdots<i_{p+1}\leq n$.
\end{thm}

\begin{prop} \label{H_to_zero}
Let $H=\left\{H_1=\cdots=H_{n-1}=0\right\}\subset \mathbb R^n$ be
a $1$-dimensional complete intersection. If $X_H$ is the
Hamiltonian vector field on $H$ then $[\mathcal
L_{X_H}(\alpha)]_{H}=[0]_{H}$ for any closed $2$-form $\alpha$.
\end{prop}

\begin{proof} Note that $X_H\rfloor dx_1\wedge\ldots\wedge
dx_n=dH_1\wedge\ldots\wedge dH_p$. This implies for $i<j$
$$X_H\rfloor dx_i\wedge dx_j=(-1)^{i+j+1}(\frac{\partial}{\partial x_{i_1}} \wedge\cdots\wedge\frac{\partial}{\partial x_{i_{n-2}}} )
\rfloor(dH_1\wedge\cdots\wedge dH_{n-1})=$$
$$=\sum_{k=1}^{n-1}(-1)^{k+i+j}(\frac{\partial}{\partial x_{i_1}} \wedge\cdots\wedge\frac{\partial}{\partial x_{i_{n-2}}} )
\rfloor(dH_{l_{1,k}}\wedge\cdots\wedge
dH_{l_{n-2,k}})dH_k=\sum_{k=1}^{n-1}f_k dH_k$$ where
$(i_1,\cdots,i_{n-2})=(1,\cdots,i-1,i+1,\cdots,j-1,j+1,\cdots,n)$
and for $k\in\{1,\cdots,n-1\}$ we take a sequence
$(l_{1,k},\cdots,l_{n-2,k})=(1,\cdots,k-1,k+1,\cdots,n-1)$.

\medskip

\noindent Thus $[X_H\rfloor dx_i\wedge
dx_j]_{H=0}=[\sum_{k=1}^{n-1}f_k dH_k]_{H}=[0]_{H}$. If
$\alpha=\sum_{i<j}g_{i,j}dx_i\wedge dx_j$ is a closed $2$-form
then $[\mathcal L_{X_H}\alpha]_H=[d(X_H\rfloor\alpha)]_H$. It
implies that $$[\mathcal
L_{X_H}\alpha]_H=\sum_{i<j}g_{i,j}[d(X_H\rfloor dx_i\wedge
dx_j)]_H+[dg_{i,j}\wedge (X_H\rfloor dx_i\wedge
dx_j)]_H=[0]_H.$$\end{proof}

The germ of a vector field tangent to $T_7$ of non trivial action on algebraic restriction of closed 2-forms to  $T_7$ may be described as a linear combination germs of vector fields: $X_0=E,\, X_1=x_3E,\, X_2=x_2E,\, X_3=x_1E,\, X_4=x_2^2E,\,
X_5=x_3^2E$ where $E$ is the Euler vector field $E= 3 x_1 \partial /\partial x_1+2 x_2 \partial /\partial x_2+2 x_3 \partial /\partial x_3$.

\begin{prop} \label{t7-infinitesimal}

The infinitesimal action of germs of quasi-homogeneous vector
fields tangent to $N\in (T_7)$ on the basis of the vector space of
algebraic restrictions of closed $2$-forms to $N$ is presented in
Table \ref{infini-t7}.

\begin{table}[h]
\begin{center}
\begin{tabular}{|l|r|r|r|r|r|r|r|}

 \hline

  $\mathcal L_{X_i} [\theta_j]$ & $[\theta_1]$   &   $[\theta_2]$ &   $[\theta_3]$ & $[\theta_4]$   & $[\theta_5]$   & $[\theta_6]$   & $[\theta_7]$ \\ \hline 

                      $X_0=E$ & $4 [\theta_1]$ & $5 [\theta_2]$ & $5 [\theta_3]$ & $7 [\theta_4]$ & $7 [\theta_5]$ & $7 [\theta_6]$ & $9 [\theta_7]$ \\ \hline 

                   $X_1=x_3E$ & $[0]$          & $7 [\theta_4]$ & $3 [\theta_6]$ & $9 [\theta_7]$ & $[0]$          & $[0]$          & $[0]$ \\      \hline

                   $X_2=x_2E$ & $[0]$          & $-3 [\theta_6]$ & $7 [\theta_5]$ & $[0]$         & $-9 [\theta_7]$ & $[0]$         & $[0]$ \\      \hline

                   $X_3=x_1E$ & $-4 [\theta_6]$ & $[0]$          & $[0]$          & $[0]$         & $[0]$           & $[0]$         & $[0]$ \\       \hline

                $X_4=x_2^2 E$ & $[0]$           & $[0]$          & $-9 [\theta_7]$ & $[0]$        & $[0]$           & $[0]$         & $[0]$ \\    \hline

                $X_5=x_3^2 E$ & $[0]$           & $9 [\theta_7]$ & $[0]$          & $[0]$         & $[0]$           & $[0]$         & $[0]$ \\  \hline

\end{tabular}
\end{center}
\caption{Infinitesimal actions on algebraic restrictions of closed
2-forms to  $T_7$. $E=3x_1 \partial /\partial x_1+ 2x_2 \partial /\partial x_2+ 2x_3 \partial /\partial x_3$}\label{infini-t7}
\end{table}

\end{prop}

\medskip

Let $\mathcal{A}=[c_1 \theta_1+c_2 \theta_2+c_3 \theta_3+c_4 \theta_4+c_5 \theta_5+c_6 \theta_6 +c_7 \theta_7]_{T_7}$
be the algebraic restriction of a symplectic form $\omega$.

The first statement of Theorem \ref{klasT7} follows from the following lemmas.

\begin{lem}
\label{t7lem0} If \;$c_1\cdot c_2\cdot c_3\ne 0$\; then the algebraic restriction
$\mathcal{A}=[\sum_{k=1}^7 c_k \theta_k]_{T_7}$
can be reduced by a symmetry of $T_7$ to an algebraic restriction $ [\theta_1+\widetilde{c}_2 \theta_2+\widetilde{c}_3 \theta_3]_{T_7}$.
\end{lem}

\begin{proof}[Proof of Lemma \ref{t7lem0}]

We use the homotopy method to prove that  $\mathcal{A}$ is diffeomorphic to $ [\theta_1+\widetilde{c}_2 \theta_2+\widetilde{c}_3 \theta_3]_{T_7}$.

Let $\mathcal{B}_t=[c_1 \theta_1+c_2 \theta_2+c_3 \theta_3+(1-t)c_4 \theta_4+(1-t)c_5 \theta_5+(1-t)c_6 \theta_6 +(1-t)c_7 \theta_7]_{T_7}$
\; for $t \in[0;1]$. Then $\mathcal{B}_0=\mathcal{A}$\; and \;$\mathcal{B}_1=[c_1 \theta_1+c_2 \theta_2+c_3 \theta_3]_{T_7}$.
 We prove that there exists a family $\Phi_t \in Symm(T_7),\;t\in [0;1]$ such that
 \begin{equation}
\label{prooft7lem01}   \Phi_t^*\mathcal{B}_t=\mathcal{B}_0,\;\Phi_0=id.
\end{equation}
Let $V_t$ be a vector field defined by $\frac{d \Phi_t}{dt}=V_t(\Phi_t)$. Then differentiating (\ref{prooft7lem01}) we obtain
  \begin{equation}
\label{prooft7lem02}   \mathcal L_{V_t} \mathcal{B}_t=c_4 \theta_4+c_5 \theta_5+c_6 \theta_6 +c_7 \theta_7.
\end{equation}
We are looking for $V_t$ in the form $V_t=\sum_{k=1}^{5}b_k(t) X_k$   where $b_k(t)$ for $k=1,\ldots,5$ are smooth functions $b_k:[0;1]\rightarrow \mathbb{R}$. Then by Proposition  \ref{t7-infinitesimal} equation (\ref{prooft7lem02}) has a form

\begin{equation}  \label{prooft7lem03}
\left[ \begin{array}{ccccc}
7c_2 & 0 & 0 & 0 & 0 \\
 0 & 7c_3 & 0 & 0 & 0 \\
3c_3 & -3c_2 & -4c_1 & 0 & 0 \\
9c_4(1-t) & -9c_5(1-t) & 0 & -9c_3 & 9c_2
\end{array} \right]
\left[ \begin{array}{c} b_1 \\ b_2 \\ b_3 \\ b_4 \\ b_5  \end{array} \right] =
\left[ \begin{array}{c} c_4 \\ c_5 \\ c_6 \\ c_7  \end{array}  \right]
\end{equation}

\noindent If \;$c_1\cdot c_2\cdot c_3\ne 0$ we can solve (\ref{prooft7lem03}) and $\Phi_t$ may be obtained as a flow of vector field $V_t$.
The family $\Phi_t$ preserves $T_7$, because $V_t$ is tangent to $T_7$ and $\Phi_t^*\mathcal{B}_t=\mathcal{A}$.
Using the homotopy arguments we have $\mathcal{A}$ diffeomorphic to $ \mathcal{B}_1=[c_1 \theta_1+c_2 \theta_2+c_3 \theta_3]_{T_7}$.
By the condition $c_1\ne 0$ we have a diffeomorphism $\Psi \in Symm(T_7)$ of the form
  \begin{equation}
\label{prooft7lem04}
\Psi:\,(x_1,x_2,x_3)\mapsto (|c_1|^{-\frac{3}{4}} x_1,|c_1|^{-\frac{1}{2}} x_2,|c_1|^{-\frac{1}{2}} x_3)
\end{equation}
and we obtain
\[ \Psi^*(\mathcal{B}_1)=[\frac{c_1}{|c_1|} \theta_1+c_2 |c_1|^{-\frac{5}{4}} \theta_2+c_3 |c_1|^{-\frac{5}{4}} \theta_3]_{T_7} =
 [\pm \theta_1+ \widetilde{c}_2 \theta_2+\widetilde{c}_3 \theta_3]_{T_7}.\]
 By the following symmetry of $T_7$: $(x_1,x_2,x_3)\mapsto (x_1, x_3, x_2)$, we have that
 $[\theta_1+ \widetilde{c}_2 \theta_2+\widetilde{c}_3 \theta_3]_{T_7}$ and $[-\theta_1+ \widetilde{c}_3 \theta_2+\widetilde{c}_2 \theta_3]_{T_7}$ are diffeomorphic.
\end{proof}

\begin{lem}
\label{t7lem1} If \;$c_2\cdot c_3= 0$\;and $c_2+ c_3 \neq 0$ then the algebraic restriction of the form
$[\sum_{k=1}^7 c_k \theta_k]_{T_7}$
can be reduced by a symmetry of $T_7$ to an algebraic restriction $ [\widetilde{c}_1 \theta_1+ \theta_2+\widetilde{c}_5 \theta_5]_{T_7}$.
\end{lem}

\begin{proof}[Proof of Lemma \ref{t7lem1}]
We use similar methods as above to prove that  if \;$c_2\cdot c_3=
0$\;and $c_2+ c_3 \neq 0$ then $\mathcal{A}$ is diffeomorphic to
$[\widetilde{c}_1 \theta_1+ \theta_2+\widetilde{c}_5
\theta_5]_{T_7}$.\par \noindent If $c_3=0$ then $c_2\neq 0$ and
$\mathcal{A}= [c_1 \theta_1+c_2 \theta_2+c_4 \theta_4+c_5
\theta_5+c_6 \theta_6 +c_7 \theta_7]_{T_7}$ Let
$\mathcal{B}_t=[c_1 \theta_1+c_2 \theta_2+(1-t)c_4 \theta_4+c_5
\theta_5+(1-t)c_6 \theta_6 +(1-t)c_7 \theta_7]_{T_7}$ \; for $t
\in[0;1]$. Then $\mathcal{B}_0=\mathcal{A}$\; and
\;$\mathcal{B}_1=[c_1 \theta_1+c_2 \theta_2+c_5 \theta_5]_{T_7}$.
 We prove that there exists a family $\Phi_t \in Symm(T_7),\;t\in [0;1]$ such that
 \begin{equation}
\label{prooft7lem11}   \Phi_t^*\mathcal{B}_t=\mathcal{B}_0,\;\Phi_0=id.
\end{equation}
Let $V_t$ be a vector field defined by $\frac{d \Phi_t}{dt}=V_t(\Phi_t)$. Then differentiating (\ref{prooft7lem11}) we obtain
  \begin{equation}
\label{prooft7lem12}   \mathcal L_{V_t} \mathcal{B}_t=c_4 \theta_4+c_6 \theta_6 +c_7 \theta_7.
\end{equation}
We are looking for $V_t$ in the form $V_t=b_1(t) X_1+b_2(t) X_2+ b_4(t)
X_4+b_5(t) X_5$   where $b_k(t)$ for $k=1,2,4,5$ are smooth functions $b_k:[0;1]\rightarrow \mathbb{R}$.
 Then by Proposition  \ref{t7-infinitesimal} equation (\ref{prooft7lem12})
has a form

\begin{equation}  \label{prooft7lem13}
\left[ \begin{array}{cccc}
7c_2 & 0 & 0 & 0 \\
0 & -3c_2 & -4c_1 & 0 \\
9c_4(1-t) & -9c_5 & 0 & 9c_2
\end{array} \right]
\left[ \begin{array}{c} b_1 \\ b_2  \\ b_4 \\ b_5  \end{array} \right] =
\left[ \begin{array}{c} c_4  \\ c_6 \\ c_7  \end{array}  \right]
\end{equation}

\noindent If \;$ c_2\ne 0$ we can solve (\ref{prooft7lem13}) and $\Phi_t$ may be obtained as a flow of vector field $V_t$.
The family $\Phi_t$ preserves $T_7$, because $V_t$ is tangent to $T_7$ and $\Phi_t^*\mathcal{B}_t=\mathcal{A}$.
Using the homotopy arguments we have that $\mathcal{A}$ is diffeomorphic to $ \mathcal{B}_1=[c_1 \theta_1+c_2 \theta_2+c_5 \theta_5]_{T_7}$.
By the condition $c_2\ne 0$ we have a diffeomorphism $\Psi \in Symm(T_7)$ of the form
  \begin{equation}
\label{prooft7lem14}
\Psi:\,(x_1,x_2,x_3)\mapsto (c_2^{-\frac{3}{5}} x_1,c_2^{-\frac{2}{5}} x_2,c_2^{-\frac{2}{5}} x_3)
\end{equation}
and we obtain
\[ \Psi^*(\mathcal{B}_1)=[c_1 c_2^{-\frac{4}{5}} \theta_1+ \theta_2+c_5 c_2^{-\frac{7}{5}} \theta_5]_{T_7} =
 [\widetilde{c}_1 \theta_1+ \theta_2+\widetilde{c}_5 \theta_5]_{T_7}.\]

\noindent If $c_2=0$ then $c_3\neq 0$ and by the  diffeomorphism  $\Theta \in Symm(T_7)$ of the form: $(x_1,x_2,x_3)\mapsto (x_1, x_3, x_2)$, we obtain
$\Theta ^*[c_1 \theta_1+c_3 \theta_3+c_4 \theta_4+c_5 \theta_5+c_6 \theta_6 +c_7 \theta_7]_{T_7} =
[-c_1 \theta_1+c_3 \theta_2+c_4 \theta_5+c_5 \theta_4-c_6 \theta_6 -c_7 \theta_7]_{T_7}$ and we may use the homotopy method now.
 \end{proof}

\begin{lem}
\label{t7lem2} If $c_2=c_3=0$, \; $c_1\neq 0$\; and \;$(c_4, c_5)\neq (0,0)$\; then the algebraic restriction of the form
$[\sum_{k=1}^7 c_k \theta_k]_{T_7}$
can be reduced by a symmetry of $T_7$ to an algebraic restriction $ [\theta_1+ \widetilde{c}_4 \theta_4+\widetilde{c}_5 \theta_5]_{T_7}$.
\end{lem}

\begin{lem}
\label{t7lem4} If \;$c_1\neq 0$\; and $c_2=c_3=c_4=c_5=0$ then the algebraic restriction of the form
$[\sum_{k=1}^7 c_k \theta_k]_{T_7}$
can be reduced by a symmetry of $T_7$ to an algebraic restriction $ [\theta_1+ \widetilde{c}_7 \theta_7]_{T_7}$.
\end{lem}


\begin{lem}
\label{t7lem3} If $c_1=c_2=c_3=0$ and \;$(c_4, c_5)\neq (0,0)$\; then the algebraic restriction of the form
$[\sum_{k=1}^7 c_k \theta_k]_{T_7}$
can be reduced by a symmetry of $T_7$ to an algebraic restriction $ [\theta_4+ \widetilde{c}_5 \theta_5+\widetilde{c}_6 \theta_6]_{T_7}$.
\end{lem}



\begin{lem}
\label{t7lem5} If $c_1=\ldots =c_5=0$ and \;$c_6\neq 0$\;  then the algebraic restriction  $\mathcal{A}=[\sum_{k=1}^7 c_k \theta_k]_{T_7}$
can be reduced by a symmetry of $T_7$ to an algebraic restriction $ [\theta_6+ \widetilde{c}_7 \theta_7]_{T_7}$.
\end{lem}


\begin{lem}
\label{t7lem6} If $c_1=\ldots =c_6=0$ and \;$c_7\neq 0$\;  then the algebraic restriction
$\mathcal{A}=[\sum_{k=1}^7 c_k \theta_k]_{T_7}$ can be reduced by a symmetry of $T_7$ to an algebraic restriction $ [\theta_7]_{T_7}$.
\end{lem}


The proofs of Lemmas \ref{t7lem2} -- \ref{t7lem6} are similar and are based on Table \ref{infini-t7}.
\smallskip

Statement $(ii)$ of Theorem \ref{klasT7} follows from conditions
in the proof of part $(i)$ and $(iii)$ follows from Theorem
\ref{geom-cond-t7} which was proved in Section \ref{t7-geom_cond}.

\smallskip

Now we prove that the parameters $c, c_1, c_2$ are moduli in the normal forms. The proofs are very similar in all cases. We consider as an example
the normal form with two parameters $[c_1\theta_1+\theta_2+c_2\theta_3]_{T_7}$. From Table \ref{infini-t7} we see that the tangent space to the orbit
of $[c_1\theta_1+\theta_2+c_2\theta_3]_{T_7}$ at $[c_1\theta_1+\theta_2+c_2\theta_3]_{T_7}$ is spanned by the linearly independent algebraic restrictions
$[4c_1\theta_1+5\theta_2+5c_2\theta_3]_{T_7}$, $[\theta_4]_{T_7},[\theta_5]_{T_7}, [\theta_6]_{T_7}, [\theta_7]_{T_7}.$ Hence the algebraic restrictions
$[\theta_1]_{T_7}$ and $[\theta_3]_{T_7}$ do not belong to it. Therefore the parameters $c_1$ and $c_2$ are independent moduli in the normal form
$[c_1\theta_1+\theta_2+c_2\theta_3]_{T_7}$.

\bigskip

\textbf{Acknowledgements.} The authors wish to express their thanks to the referee for many valuable suggestions.


\bibliographystyle{amsalpha}

\begin{thebibliography}{AAAA}

\bibitem [A1] {Ar1} V. I. Arnold, \emph{Simple singularities of curves}, Proc. Steklov Inst. Math. 1999, no. 3 (226),
20-28.

\bibitem [A2] {Ar2} V. I. Arnold, \emph{ First step of local symplectic algebra},
Differential topology, infinite-dimensional Lie algebras, and
applications. D. B. Fuchs' 60th anniversary collection.
Providence, RI: American Mathematical Society. Transl., Ser. 2,
Am. Math. Soc. 194(44), 1999,1-8.

\bibitem [AG] {ArGi} V. I. Arnold, A. B. Givental \emph{Symplectic geometry}, in Dynamical systems, IV, 1-138, Encyclopedia of Matematical Sciences, vol. 4,
Springer, Berlin, 2001.

\bibitem [AGLV] {ArGorLyVa} V. I. Arnold, V.V. Goryunov, O.V. Lyashko, V.A. Vasil'ev \emph{Singualrity Theory II, Classification and Applications} in Dynamical systems VIII, Encyclopedia Math. Sci 39, Springer, Berlin, 1993,33-34.

\bibitem [AVG] {AVG} V. I. Arnold, S. M. Gusein-Zade, A. N. Varchenko,
\emph{ Singularities of Differentiable Maps}, Vol. 1, Birhauser,
Boston, 1985.

\bibitem [D] {D} W. Domitrz,
\emph{Local symplectic algebra of quasi-homogeneous curves}, Fundamentae Mathematicae 204 (2009), 57-86.

\bibitem [DJZ1] {DJZ1} W. Domitrz, S. Janeczko, M. Zhitomirskii,
\emph{Relative Poincare lemma, contractibility, quasi-homogeneity
and vector fields tangent to a singular variety}, Ill. J. Math.
48, No.3 (2004), 803-835.

\bibitem [DJZ2] {DJZ2} W. Domitrz, S. Janeczko, M. Zhitomirskii,
\emph{Symplectic singularities of varietes: the method of
algebraic restrictions}, J. reine und angewandte Math. 618 (2008),
197-235.

\bibitem[DR] {DR} W. Domitrz, J. H. Rieger, \emph{Volume preserving subgroups of
$\mathcal A$ and $\mathcal K$ and singularities in unimodular
geometry}, Mathematische Annalen
   345(2009), 783–-817.

\bibitem[DT] {DT} W. Domitrz, Z. Trebska, \emph{Symplectic $S_{\mu}$ singularities}, preprint, arXiv:1101.5176.

\bibitem [G] {G} M. Giusti, \emph{Classification des singularit\'es isol\'ees d'intersections compl\`etes simples},
C. R. Acad. Sci., Paris, Sér. A 284 (1977), 167-170.

\bibitem [IJ1] {IJ1} G. Ishikawa, S. Janeczko, \emph{ Symplectic bifurcations of plane curves and isotropic liftings},
Q. J. Math. \textbf{54}, No.1 (2003), 73-102.

\bibitem [IJ2] {IJ2} G. Ishikawa, S. Janeczko, \emph{ Symplectic singularities of isotropic mappings},
Geometric singularity theory, Banach Center Publications
\textbf{65} (2004), 85-106.

\bibitem [K] {K} P. A. Kolgushkin, \emph{Classification of simple multigerms of curves
in a space endowed with a symplectic structure}, St. Petersburg
Math. J. \textbf{15} (2004), no. 1, 103-126.

\bibitem [Wa] {Wa} J. M. Wahl, \emph{Derivations, automorphisms and deformations of quasi-homogeneous
singularities}, Singularities, Summer Inst., Arcata/Calif. 1981,
Proc. Symp. Pure Math. 40, Part 2, 613-624 (1983).

\bibitem [W] {W} C. T. C. Wall,
\emph{Singular points of plane curves}, London Mathematical
Society Student Texts, 63, Cambridge University Press, Cambridge,
2004.

\bibitem [Z] {Zh} M. Zhitomirskii, {\em Relative Darboux
theorem for singular manifolds and local contact algebra}, Can. J.
Math. \textbf{57}, No.6 (2005), 1314-1340.
\end{thebibliography}

\end{document}